\documentclass[10pt,a4paper]{article}
\usepackage[latin1]{inputenc}
\usepackage{ucs}
\usepackage{amsmath}
\usepackage{amsfonts}
\usepackage{amssymb}
\usepackage{upgreek}
\usepackage{graphicx}
\usepackage{pgfplots}
\usepackage{bm}
\usepackage[bookmarksnumbered,colorlinks, linkcolor=black, citecolor=black]{hyperref}

\usepackage{pdfsync}
\usepackage{srcltx}
\author{age}

\newcommand{\red}[1]{\textcolor{black}{#1}} 
\newcommand\vm[1]{\bm{\mathrm{#1}}} 
\newcommand{\norm}[1]{\lVert#1\rVert}

\newcommand{\Eref}[1]{(\ref{#1})}

\newcommand{\xx}{\mathbf{x}}
\newcommand{\qq}{\mathbf{q}}

\newcommand{\bb}{\mathbf{B}}
\newcommand{\KK}{\mathbf{K}}
\newcommand{\DD}{\mathbf{D}}

\pgfplotscreateplotcyclelist{ageplot}{%
{blue,mark=o, mark size={3}},
{red,mark=square},
{black,mark=star},
{brown,mark=triangle},
{green,mark=diamond},}

\begin{document}


\begin{center}
\large \red{Efficient recovery-based error estimation} for the smoothed finite element method for smooth and singular linear elasticity 
\end{center}

\begin{center}{Octavio A.~González-Estrada$^{1}$, Sundararajan Natarajan$^{1}$,  Juan~Jos\'{e}~R\'{o}denas$^{2}$, Hung Nguyen-Xuan$^{3}$, St\'{e}phane~P.A.~Bordas$^{1}$\footnote{Correspondence to: stephane.bordas@northwestern.edu} } 
\end{center}

\begin{center}\small{
$^{1}$Institute of Mechanics \& Advanced Materials, School of Engineering, Cardiff University, Queen's Building, The Parade, Cardiff CF24 3AA Wales, UK.\\ 
$^{2}$Centro de Investigación de Tecnología de Vehículos(CITV),
\\Universitat Polit\`{e}cnica de Val\`{e}ncia, E-46022-Valencia, Spain.}\\
$^{3}$ Dept. Mechanics, Faculty of Mathematics and Computer Science in the University of Science at Vietnam National University in Ho Chi Minh (HCM) City. \\
\end{center}

\begin{abstract}

An error control technique aimed to assess the quality of smoothed finite element approximations is presented in this paper. Finite element techniques based on strain smoothing appeared in 2007 were shown to provide significant advantages compared to conventional finite element approximations. In particular, a widely cited strength of such methods is improved accuracy for the same computational cost. Yet, few attempts have been made to directly assess the quality of the results obtained during the simulation by evaluating an estimate of the discretization error. \red{Here we propose a recovery type error estimator based on an enhanced recovery technique. The salient features of the recovery are: enforcement of local equilibrium and, for singular problems a ``smooth+singular'' decomposition of the recovered stress. We evaluate the proposed estimator on a number of test cases from linear elastic structural mechanics and obtain precise error estimations whose effectivities, both at local and global levels, are improved compared to recovery procedures not implementing these features. 
}
\end{abstract}

\begin{footnotesize}
\end{footnotesize}

\section{Introduction}

The smoothed finite element for mechanics problems was introduced in 2006 by Liu \textit{et al.,}~\cite{liudai2006}. The main idea of the method is to relax the compatibility condition at the element level by replacing the standard compatible strain by its smoothed counterpart. The smoothing operation can be performed over domains of various shapes which can be obtained by dividing the computational domain into non-overlapping smoothing domains. These domains can be obtained by subdividing the elements (cell-based smoothing) as in \cite{liudai2006,liunguyen2007,nguyen-xuanbordas2008,bordasnatarajan2010,zhangliu2008}, or using edge \cite{nguyen-thoiliu2009,liunguyen-thoi2009b} or node-based geometrical information \cite{liunguyen-thoi2009}. Each method has several advantages and drawbacks, summarized in, e.g. \cite{liu2010,liunguyen-xuan2010}, but the strongest motivation for smoothed finite elements is certainly revealed in its enhancement of low order simplex elements (e.g. linear triangles and tetrahedral), alleviating overstiffness, locking and improving their accuracy significantly \cite{nguyenliu2007, hungbordas2009, nguyen-thoiliu2009}.

The applications of strain smoothing in finite elements are wide. Since the introduction of the smoothed finite element method (SFEM), the convergence, the stability, the accuracy and the computational complexity of this method was studied in~\cite{liunguyen2007,nguyen-xuanbordas2008} and the method was extended to treat various problems in solid mechanics such as plates~\cite{nguyen-xuanrabczuk2008}, shells~\cite{nguyenrabczuk2008} and nearly incompressible elasticity~ \cite{hungbordas2009}. Recently, Bordas \textit{et al.}~\cite{bordasrabczuk2010, bordasnatarajan2011} combined strain smoothing with the XFEM to obtain the Smoothed eXtended Finite Element Method to solve problems with strong and weak discontinuities in 2D continuum.

In this paper, we focus on the cell-based smoothing, a review of which is provided in \cite{bordasnatarajan2010,hungbordas2009,liunguyen-thoi2009a} along with applications to plates, shells, three dimensional continuum and a coupling with the extended finite element method with applications to linear elastic fracture (continuum, plate). 


The development of new numerical techniques based on the finite element method (e.g. GFEM, XFEM, ...) aims at obtaining more accurate solutions for engineering problems. Despite the accuracy gained with the new techniques, numerical errors, especially the discretization error, are always present and have to be evaluated. Consolidated accuracy assessment techniques developed in the FE framework are commonly adapted to the framework of these new techniques \cite{stroubouliszhang2006,bordasduflot2007,xiaokarihaloo2004,rodenasgonzalez2008,panetierladeveze2010}. As in any numerical method, the smoothed FEM approximation introduces an error that needs to be controlled to guarantee the quality of the numerical simulations. Although an adaptive node-based smoothed FEM has been developed in \cite{nguyen-thoiliu2011}, a rather simple error estimator using a recovery procedure which initially is only valid for NS-FEM is used to guide the adaptive process. The technique evaluates a first-order recovered strain field interpolating the nodal values by means of the linear FEM shape functions. 

The urge for quality assessment tools for smoothed FEM approximations, the promising results in \cite{nguyen-thoiliu2011} obtained with a rather simple technique and the experience of the authors in the development of high quality recovery-based error estimators in the FEM and XFEM contexts, motivates this paper where we estimate, a posteriori, the approximation error committed by cell-based smoothed finite elements. The method used for a posteriori error estimation relies on the Zienkiewicz and Zhu error estimator \cite{zienkiewiczzhu1987} commonly used in FEM, together with a recovery technique recently developed by the authors, specially tailored to the analysis of enriched approximations containing a smooth and a singular part and which locally enforces the fulfilment of equilibrium equations. The technique known as SPR-CX \cite{rodenasgonzalez2008, rodenasgonzalez2010} was shown to lead to very good effectivity indices in FEM and XFEM. 

The paper is organised as follows: In Section \ref{sec:ProbStatement}, the boundary value problem of linear elasticity is briefly introduced and the approximate solution using the smoothed finite element method is presented. In Section \ref{sec:Error}, we discuss basic concepts related to error estimation, in especial recovery-based techniques. Section  \ref{section:sprcx} is devoted to the proposed enhanced recovery technique and its application to SFEM approximations. Numerical examples are presented in Section \ref{sec:Results} and the main concluding remarks in Section \ref{sec:Conclusions}.

\section{Problem statement and SFEM solution}
\label{sec:ProbStatement}

\subsection{ Problem statement}
 
Let us consider the 2D linear elasticity problem. The unknown displacement field $\vm{u}$, taking values in $\Omega \subset \mathbb{R}^{2}$, is the solution of the boundary value problem given by 
\begin{align}
  -\nabla \cdot \vm{\upsigma} \left(\vm{u}\right) &= \vm{b}  	&&  {\rm in }\; \Omega 
   \label{Eq:IntEq} \\
   \vm{\upsigma} \left(\vm{u} \right)\cdot \vm{n} &= \vm{t} 	&&  {\rm on }\; \Gamma _{N}  		\label{Eq:Neumann}\\
   \vm{u}                                         &= \vm{0}	&&  {\rm on }\; \Gamma _{D} \label{Eq:Dirichlet}
\end{align}
 
where $\Gamma _{N}$ and $\Gamma _{D}$ denote the Neumann and Dirichlet boundaries with $\partial \Omega = \Gamma_N \cup \Gamma_D $ and $\Gamma_N \cap \Gamma_D  =\emptyset$. The Dirichlet boundary condition in \Eref{Eq:Dirichlet} is assumed to be homogeneous for the sake of simplicity. 
 
The weak form of the problem reads: Find $\vm{u} \in V$ such that 
\begin{equation} \label{Eq:WeakForm} 
\forall \vm{v} \in V \qquad a(\vm{u},\vm{v}) = l(\vm{v}),
\end{equation}  
 
where ${V}$ is the standard test space for the elasticity problem such that $V = \{\vm{v} \;|\; \vm{v} \in  [H^1(\Omega)]^2 , \vm{v}|_{\Gamma_D}(\vm{x}) = \vm{0} \}$, and 
\begin{align}
a(\vm{u},\vm{v})  & := \int _{\Omega}  \vm{\varepsilon}(\vm{u})^T  \vm{D}  \vm{\varepsilon}(\vm{v}) d \Omega =
\int _{\Omega}  \vm{\sigma}(\vm{u})^T  \vm{D}^{-1}  \vm{\sigma}(\vm{v}) d \Omega \\
l(\vm{v})&:=\int _{\Omega} \vm{b}^T  \vm{v}d \Omega + \int _{\Gamma_N} \vm{t}^T \vm{v}d \Gamma,
\end{align}
 
where $\vm{D} $ is the elasticity matrix of the constitutive relation $\vm{\sigma}= \vm{D} \vm{\varepsilon} $ , $ \vm{\sigma}$ and $\vm{\varepsilon}$ denote the stress and strain operators.

\subsubsection{Singular problem}

\begin{figure}[!htb]
	\centering
	\includegraphics{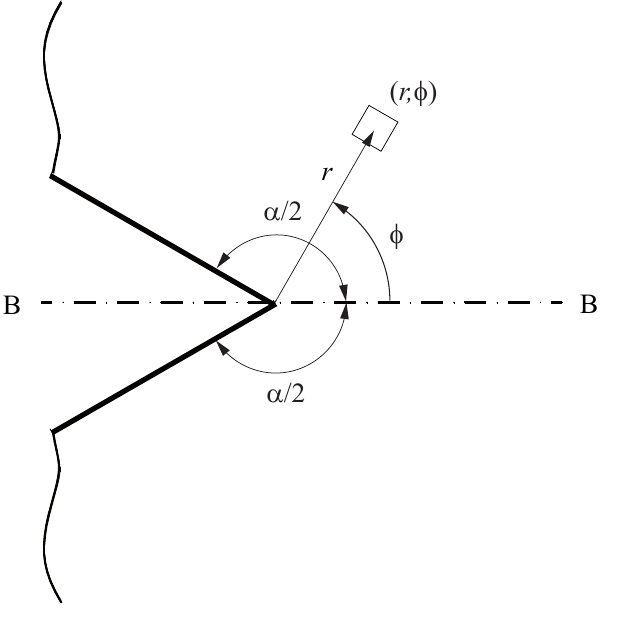} 
	\caption{Sharp reentrant corner in an infinite half-space.}
	\label{fig:Vnotch}
\end{figure}

Figure \ref{fig:Vnotch} shows a portion of an elastic body with a reentrant corner (or V-notch), subjected to tractions on remote boundaries. No body loads are applied. For this kind of problem, the stress field exhibits a singular behaviour at the notch vertex.

The analytical solution corresponding to the stress distribution in the vicinity of the singular point is a linear combination of the singular and the non-singular terms. It is often claimed that the term with a highest order of singularity dominates over the other terms in a zone surrounding the singular point sufficiently closely. The analytical solution to this singular elastic problem in the vicinity of the singular point was first given by \cite{williams1952} and is described, for example, in \cite{szabobabuska1991,barber2010}. Here, we reproduce those expressions for completeness. In accordance with the polar coordinate system of Figure \ref{fig:Vnotch}, the displacement and the stress fields at points sufficiently close to the corner can be described as:
\begin{align} 
\vm{u}(r,\phi ) &= 
K_{{\rm I}}  r^{\lambda _{{\rm I}} } \vm{\Psi} _{{\rm I}} (\lambda _{{\rm I}} ,\phi ) + 
K_{{\rm II}}  r^{\lambda _{{\rm II}} }  \vm{\Psi} _{{\rm II}} (\lambda _{{\rm II}} ,\phi ) \label{Eq:VnotchDisp} \\ 
\vm{\sigma} (r,\phi ) &= 
K_{{\rm I}}  \lambda _{{\rm I}}  r^{\lambda _{{\rm I}} -1}  \vm{\Phi} _{{\rm I}} (\lambda _{{\rm I}} ,\phi ) + 
K_{{\rm II}}  \lambda _{{\rm II}} r^{\lambda _{{\rm II}} -1} \vm{\Phi} _{{\rm II}} (\lambda _{{\rm II}} ,\phi ) \label{Eq:VnotchStress}
\end{align} 

where $r$ is the radial distance to the corner, $\lambda _m$ (with $m = {\rm I,II}$) are the eigenvalues that determine the order of the singularity, $ \vm{\Psi} _m$ and $ \vm{\Phi}_m $ are a set of trigonometric functions that depend on the angular position $\phi$, and $K_m $ are the so-called generalised stress intensity factors (GSIFs). The generalized stress intensity factor (GSIF) is a multiplicative constant that depends on the loading of the problem and linearly determines the intensity of the displacement and the stress fields in the vicinity of the singular point. Hence, the eigenvalues $\lambda$ and the GSIFs $K$ define the singular field.

The eigenvalue $\lambda$ is easily known because it depends only on the corner angle $\alpha$, and can be obtained as the smallest positive root of the following characteristic equations \cite{williams1952}:
\begin{align}  
\sin  \lambda _{{\rm I}} \alpha  + \lambda _{{\rm I}} \sin \alpha  &= 0 \label{Eq:VnotchEigen1} \\ 
\sin  \lambda _{{\rm II}} \alpha - \lambda _{{\rm II}} \sin \alpha &= 0 \label{Eq:VnotchEigen2}
\end{align} 

The smallest positive root yields the highest order of singularity and determines the term that dominates the elastic fields given by \Eref{Eq:VnotchDisp} in the vicinity of the notch vertex. Strictly speaking, \Eref{Eq:VnotchEigen1} corresponds to the symmetric part of the elastic fields with respect to $\phi = 0$ (i.e. the bisector line BB in Figure~\ref{fig:Vnotch}) and \Eref{Eq:VnotchEigen2} to the antisymmetric solution. These solutions are also called mode I and mode II solutions, respectively. The trigonometric functions for the mode I displacement and stress fields in (\ref{Eq:VnotchDisp}, \ref{Eq:VnotchStress}) are given by \cite{szabobabuska1991}:
\begin{align}
	\vm{\Psi}_{\rm I}(\lambda_{\rm I},\phi) &= 
	\begin{Bmatrix}
	  \vm{\Psi}_{{\rm I},x}(\lambda_{\rm I},\phi)\\
	  \vm{\Psi}_{{\rm I},y}(\lambda_{\rm I},\phi)
	\end{Bmatrix} \nonumber\\ &=
	\frac{1}{2 \mu}
	\begin{Bmatrix}
	  (\kappa -Q (\lambda_{\rm I} +1))\cos \lambda_{\rm I} \phi -\lambda_{\rm I} \cos(\lambda_{\rm I} -2)\phi) \\
	  (\kappa +Q (\lambda_{\rm I} +1))\sin \lambda_{\rm I} \phi +\lambda_{\rm I} \sin(\lambda_{\rm I} -2)\phi
 	\end{Bmatrix} \label{Eq:TrigoPsiI}\\
	\vm{\Phi}_{\rm I}(\lambda_{\rm I},\phi) &=
	\begin{Bmatrix}
	  \vm{\Phi}_{{\rm I},xx}(\lambda_{\rm I},\phi)\\
	  \vm{\Phi}_{{\rm I},yy}(\lambda_{\rm I},\phi)\\
	  \vm{\Phi}_{{\rm I},xy}(\lambda_{\rm I},\phi)
	\end{Bmatrix} \nonumber\\ &=
	\begin{Bmatrix}
	  (2-Q(\lambda_{\rm I}+1))\cos(\lambda_{\rm I}-1)\phi - (\lambda_{\rm I}-1)\cos(\lambda_{\rm I}-3)\phi\\
	  (2+Q(\lambda_{\rm I}+1))\cos(\lambda_{\rm I}-1)\phi + (\lambda_{\rm I}-1)\cos(\lambda_{\rm I}-3)\phi \\
	  Q(\lambda_{\rm I}+1)\sin(\lambda_{\rm I}-1)\phi + (\lambda_{\rm I}-1)\sin(\lambda_{\rm I}-3)\phi
	\end{Bmatrix} \label{Eq:TrigoPhiI}
\end{align}

where $\mu$ is the shear modulus and $\kappa$ is the Kolosov constant, defined as functions of $E$ (Young's modulus) and $\nu$ (Poisson's coefficient) according to the following expressions:
\begin{equation*}
\mu =\frac{E}{2\left(1+\nu \right)}, \qquad 
\kappa =
  \begin{cases}
    3-4 \nu                & \textrm{plane strain} \\
   \dfrac{3-\nu}{1+\nu}    & \textrm{plane stress}
  \end{cases}
\end{equation*}

In the same way, for mode II we have:
\begin{align}
	\vm{\Psi}_{\rm II}(\lambda_{\rm II},\phi) &= 
	\begin{Bmatrix}
	  \vm{\Psi}_{{\rm II},x}(\lambda_{\rm II},\phi)\\
	  \vm{\Psi}_{{\rm II},y}(\lambda_{\rm II},\phi)
	\end{Bmatrix} \nonumber\\ &=
	\frac{1}{2 \mu}
	\begin{Bmatrix}
	  (\kappa -Q (\lambda_{\rm II} +1))\sin \lambda_{\rm II} \phi -\lambda_{\rm II} \sin(\lambda_{\rm II} -2)\phi\\
	  -(\kappa +Q (\lambda_{\rm II} +1))\cos \lambda_{\rm II} \phi -\lambda_{\rm II} \cos(\lambda_{\rm II} -2)\phi
	\end{Bmatrix} \\
	\vm{\Phi}_{\rm II}(\lambda_{\rm II},\phi) &=
	\begin{Bmatrix}
	  \vm{\Phi}_{{\rm II},xx}(\lambda_{\rm II},\phi)\\
	  \vm{\Phi}_{{\rm II},yy}(\lambda_{\rm II},\phi)\\
	  \vm{\Phi}_{{\rm II},xy}(\lambda_{\rm II},\phi)
	\end{Bmatrix} \nonumber\\ &=
	\begin{Bmatrix}
	  (2-Q(\lambda_{\rm II}+1))\sin(\lambda_{\rm II}-1)\phi - (\lambda_{\rm II}-1)\sin(\lambda_{\rm II}-3)\phi\\
	  (2+Q(\lambda_{\rm II}+1))\sin(\lambda_{\rm II}-1)\phi + (\lambda_{\rm II}-1)\sin(\lambda_{\rm II}-3)\phi\\
	  -Q(\lambda_{\rm II}+1)\cos(\lambda_{\rm II}-1)\phi + (\lambda_{\rm II}-1)\cos(\lambda_{\rm II}-3)\phi
	\end{Bmatrix} 
\end{align}

Note that the components of the displacement and the stress fields are expressed in Cartesian coordinates. In addition, $Q$ is a constant for a given notch angle:
\begin{equation} \label{Eq:QConstant}
Q_{\rm I}=-\dfrac{\cos \left((\lambda_{\rm I} -1)\dfrac{\alpha }{2} \right)}{\cos \left((\lambda_{\rm I} +1)\dfrac{\alpha }{2} \right)} \quad,\quad 
Q_{\rm II}=-\dfrac{\sin \left((\lambda_{\rm II} -1)\dfrac{\alpha }{2} \right)}{\sin \left((\lambda_{\rm II} +1)\dfrac{\alpha }{2} \right)} \qquad
\end{equation}

\subsection{FEM solution with strain smoothing}

\subsubsection{Finite element formulation}

Let $\vm{u}^{h}$ be a finite element approximation to  $\vm{u}$. The solution lies in a functional space $V^{h} \subset V$ associated with a mesh of isoparametric finite elements of characteristic size $h$, and it is such that 
\begin{equation}  
\forall \vm{v}^h \in V^{h} \qquad a(\vm{u}^{h},\vm{v}^h) = l(\vm{v}^h)
\end{equation}  

Using a variational formulation of the problem in (\ref{Eq:IntEq}--\ref{Eq:Dirichlet}) and a finite element approximation $\vm{u}^{h} =\vm{N}\vm{u}^{e} $, where $\vm{N}$ denotes the shape functions of order $p$, we obtain a system of linear equations to solve the displacements at nodes  $\vm{u}^{e}$: 
\begin{equation}    
\vm{K} \vm{U} = \vm{f}  
\end{equation}  
 
\noindent where $\vm{K}$ is the stiffness matrix, $\vm{U}$ is the vector of nodal displacements and $\vm{f} $ is the load vector. 

\subsubsection{Strain smoothing in FEM}
Inspired by the work of Chen \textit{et al.,} in 2001~\cite{chenwu2001} on stabilized conforming nodal integration (SCNI), Liu \textit{et al.,}~\cite{liudai2006} introduced the smoothed finite element method. The idea behind SCNI/SFEM is to write a strain measure as a spatial average of the standard strain field. To do so, the elements are divided into smoothing cells over which the strain is smoothed. By the divergence theorem, integration over the element is transformed to contour integration around the boundaries of the subcell. A particular element can have a certain number of smoothing cells and depending on that number, the formulation offers a range of different properties~\cite{liudai2006,liunguyen2007,nguyen-xuanbordas2008}. Interested readers are referred to the literature~\cite{liunguyen2007,nguyen-xuanbordas2008,hungbordas2009} for a detailed description of the method, its variants~\cite{liunguyen-thoi2009b,liunguyen-thoi2009a,liunguyen-thoi2009} and its convergence properties~\cite{nguyen-xuanbordas2008,nguyenliu2007}. The strain field $\tilde {\varepsilon }_{ij}^h$, used to evaluate the stiffness matrix is computed by a weighted average of the standard strain field $ {\varepsilon }_{ij}^h$. At a point $\xx_C$ in an element $\Omega^h$,
\begin{equation}
\tilde {\varepsilon }_{ij}^h (\xx_{C} )=\int_{\Omega ^h}
{\varepsilon _{ij}^h (\xx)\Phi (\xx-\xx_{C} )d\xx }
\label{eqn:epsilonvar}
\end{equation}
\noindent where $\Phi $ is a smoothing function that generally
satisfies the following properties~\cite{yoomoran2004}
\begin{equation}
\Phi \geq 0 ~~~~ \text{and} ~~~~ \int_{\Omega ^h} {\Phi(\xx) d\xx
}=1 \label{eqn:propersmoothfunc}
\end{equation}

One possible choice of $\Phi$ is given by:
\begin{equation}
\Phi = \frac{1}{A_C} ~~~~ \text{in} ~~~~ \Omega_C \quad \text{and}
\quad \Phi = 0 \quad \text{elsewhere}
 \label{eqn:propersmoothfunc2}
\end{equation}

\noindent where $A_C$ is the area of the subcell. To use~\Eref{eqn:epsilonvar}, the subcell containing point $\xx_C$
must first be located in order to compute the correct value of the
weight function $\Phi$.

The discretized strain field is obtained through the smoothed discretized gradient operator or the smoothed strain displacement operator, $\tilde{\bb}_C$, defined by
\begin{equation}
\tilde{\boldsymbol{\varepsilon}}^h(\xx_C)=\tilde{\bb}_C(\xx_C)\qq
\end{equation}

where $\qq$ are the unknown displacements coefficients defined at the nodes of the finite element, as usual. The smoothed element stiffness matrix for element $e$ is computed by the \emph{sum of the contributions of the subcells}\footnote{The subcells $\Omega_C$ form a partition of the element $\Omega^h$.}

\begin{equation}
\tilde {\KK}^e =\sum\limits_{C=1}^{nc}\int_{\Omega_C} {\tilde
{\bb}_C^\textup{T} \DD\tilde {\bb}_C } d\Omega =\sum\limits_{C=1}^{nc}
{\tilde {\bb}_C^\textup{T} \DD\tilde {\bb}_C } \int_{\Omega_C}d\Omega
=\sum\limits_{C=1}^{nc}{\tilde {\bb}^\textup{T}_C \DD\tilde {\bb}_CA_C }
\label{eqn:stiffnessvar}
\end{equation}

where $nc$ is the number of the smoothing cells of the element. The strain displacement matrix $\tilde \bb_{C}$ is constant over each $\Omega_C$ and is of the following form

\begin{equation}
\tilde \bb_{C}=\left[ {\begin{array}{*{20}c} \tilde {\bb}_{C1} &
\tilde {\bb}_{C2} & \tilde {\bb}_{C3} & \tilde {\bb}_{C4}
\end{array}} \right] \label{eq:Btilde3dbis}
\end{equation}

where for all shape functions $I \in \{1,\dots,4\}$, the $3 \times 2$ submatrix $\tilde {\bb}_{CI}$ represents the contribution to the strain displacement matrix associated with shape function $I$ and cell $C$ and  writes

\begin{equation}
\renewcommand\arraystretch{2}
\forall I\in\{1,2,\dots,4\},\forall C \in \{1,2,\dots nc \}\tilde
{\bb}_{CI} =
\int_{S_C} \left[
\begin{array}{*{20}c}
 { n_x} & 0 \\
 0 & { n_y} \\
 { n_y} & { n_x} \\
\end{array}  \right](\xx) N_I(\xx) dS
\label{eqn:Btilde3d}
\end{equation}

or, since~\Eref{eqn:Btilde3d} is computed on the boundary of $\Omega_C$ and one Gau\ss{} point is sufficient for an exact integration (in the case of a bilinear approximation):

\begin{equation}
\renewcommand\arraystretch{2}
\tilde {\bb}_{CI}(\xx_C)=\frac{1}{A_C} \sum\limits_{b=1}^{nb} {\left(
{{\begin{array}{*{20}c}
 {N_I\left(\xx_b^G\right) n_x} & 0 \\
 0 & {N_I\left(\xx_b^G\right) n_y} \\
 {N_I\left(\xx_b^G\right) n_y} & {N_I\left(\xx_b^G\right) n_x} \\
\end{array} }} \right)} l_b^C
\label{eqn:equation3.14}
\end{equation}

where $n_b$ is number of edges of the subcell, $(n_x,n_y)$ is the outward normal to the smoothing cell, $\Omega_C$, $\xx_b^G$ and ${l}_b^C$ are the center point (Gau\ss{} point) and the length of $\Gamma_b^C$, respectively.

\section{Error estimation in the energy norm}
\label{sec:Error}

The discretization error in the standard finite element approximation is defined as the difference between the exact solution $\vm{u}$ and the finite element solution $\vm{u}^{h}$: $ \vm{e} = \vm{u} -\vm{u}^{h} $. Since the exact solution is in practice unknown, in general, the exact error can only be estimated. To obtain an estimation of $\vm{e}$, norms that allow a better global interpretation of the error are normally used. The Zienkiewicz-Zhu \cite{zienkiewiczzhu1987} error estimator  defined as
 \begin{equation} \label{Eq:ZZestimator} 
\lVert\vm{e}\lVert ^{2} \approx \lVert \vm{e}_{es}\lVert^{2}=\int _{\Omega}\left( \vm{\sigma}^*- \vm{\sigma}^h \right)^{T} \vm{D}^{-1} \left(\vm{\sigma}^*- \vm{\sigma}^h \right)d\Omega   
\end{equation} 
 
relies on the recovery of an improved stress field $\vm{\sigma} ^{*} $, which is supposed to be more accurate than the FE solution $\vm{\sigma} ^{h} $, to obtain an estimation of the error in energy norm $\lVert \vm{e}_{es}\lVert$. The domain $\Omega$ could refer to the full domain of the problem or a local subdomain (element).
   
The recovered stress field $\vm{\sigma}^*$ is usually interpolated in each element using the shape functions $\vm{N}$ of the underlying FE approximation and the values of the recovered stress field calculated at the nodes ${\vm{\sigma}}^*$, given by:
\begin{equation} \label{Eq:rec-stress-interpolation}  
\vm{\sigma}^*(\vm{x}) = \sum _{I=1}^{n_e} N_{I}(\vm{x}) {\vm{\sigma}}^*_I(\vm{x}_I),
\end{equation}

where $n_e$ is the number of nodes in the element under consideration and ${\vm{\sigma}}^*_I(\vm{x}_I)$ are the stresses provided by a recovery technique at node $I$. The superconvergent patch recovery technique (SPR) proposed by Zienkiewicz and Zhu \cite{zienkiewiczzhu1992} is commonly used to evaluate the components ($j=xx,yy,xy$) of ${\vm{\sigma}}^*_I$ using a polynomial expansion,  ${\sigma}_{I,j} = \vm{p}\vm{a}$. This expansion is defined over a set of contiguous elements connected to node $I$ called \textit{patch}, where $\vm{p}$ is the polynomial basis and $\vm{a}$ are the unknown coefficients obtained using a least squares fitting to the values of the FE stresses evaluated at integration points in the patch, being $p$, normally, of the same order as the interpolation of displacements. The ZZ error estimator is asymptotically exact (i.e. the approximate error converges to the exact error as the mesh size goes to zero) if the recovered solution used in the error estimation converges at a higher rate than the finite element solution \cite{zienkiewiczzhu1992, zienkiewiczzhu1992a}. As it can be seen in (\ref{Eq:ZZestimator}), the accuracy of the error estimate is closely related to the quality of the recovered field. For this reason, several techniques have been developed aiming to improve the quality of $\vm{\sigma}^*$. The authors have proposed different techniques mostly for the FEM/XFEM context as the extended moving least squares recovery (XMLS) and the extended global recovery techniques proposed by Duflot and Bordas in a series of papers \cite{bordasduflot2007, duflotbordas2008, bordasduflot2008}, the SPR-C and the SPR-CX by R\'odenas \emph{et al.} \cite{rodenastur2007, rodenasgonzalez2008}, which were used later as the basis for the development of recovery-based error bounding techniques \cite{rodenasgonzalez2010, diezrodenas2007}. 

The next section presents the SPR-CX technique which improves the recovered field by enforcing equilibrium and effectively dealing with singular fields. 

\red{\noindent \textbf{Remark:} In mathematics is common to consider that one can only speak about an \emph{error estimator} if sharp or at least approximated upper - and desired - also lower error bounds can be proven, reserving the word \emph{indicator} when the technique does not necessarily bound the error. However, this terminology is not general and many other authors, usually from the engineering community, use the term \emph{error estimator} even when the technique is not able to provide error bounds. This is the case for example in \cite{zienkiewiczzhu1987,wibergabdulwahab1993,blackerbelytschko1994} and also our case.}

\section{SPR-CX recovery technique} \label{section:sprcx}

The SPR-CX recovery technique first introduced by R\'odenas \emph{et al.}  \cite{rodenasgonzalez2008} is an enhancement of the superconvergent patch recovery technique in \cite{zienkiewiczzhu1992}, which incorporates the ideas proposed in \cite{rodenastur2007} to guarantee the exact satisfaction of the equilibrium locally on patches. In \cite{rodenasgonzalez2008, rodenasgonzalez2010} a set of key ideas are proposed to modify the standard SPR allowing its use with singular problems. 

\begin{figure}[!htb]
	\centering
	\includegraphics{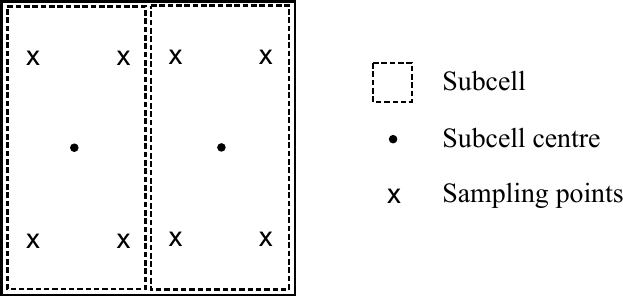}
	\caption{Distribution of stress sampling points at each subcell in a 2-subcell quadrilateral element used in the stress projection.}
	\label{fig:SC2GPmapping}
\end{figure}

The recovered stresses $\vm{\sigma}^*$ are directly evaluated at a sampling point (e.g an integration point) $\vm{x}$ through the use of a partition of unity procedure, properly weighting the stress interpolation polynomials obtained from the different patches formed at the vertex nodes of the element containing $\vm{x}$:
\begin{equation} \label{Eq:conjoint_polynomials} 
\vm{\sigma}^*(\vm{x}) = \sum _{I=1}^{n_v} N_{I}(\vm{x}) \vm{\sigma}^*_I(\vm{x}),
\end{equation}

where $N_I$ are the shape functions associated to the vertex nodes $n_v$. To obtain the nodal values $\vm{\sigma}_I$, we solve a least squares approximation of the stresses evaluated at a set of sampling points distributed within the domain of the patch of node $I$ (elements connected to $I$). In FEM,  such points usually correspond to the integrations points used in the finite element approximation. In SFEM, we map the constant strains at each subcell to a $2\times2$ Gauss quadrature distribution in the subcell used as sampling points. This way we have a sufficient number of points at each patch to solve the linear system of the least squares approximation, see Figure \ref{fig:SC2GPmapping}. Note that as in the other versions of the SFEM (NS-FEM, ES-FEM) the elements are also subdivided into subcells, a similar approach can be used to perform the mapping of the stresses to sampling points. Therefore, the proposed error estimation technique can be used with all SFEM implementations. 

One major modification of the original SPR technique is the introduction of a splitting procedure to perform the recovery. For singular problems the exact stress field $\vm{\sigma}$  is decomposed into two stress fields, a smooth field $\vm{\sigma}_{smo}$ and a singular field $\vm{\sigma}_{sing}$:
\begin{equation} \label{Eq:splitting} 
\vm{\sigma} =   \vm{\sigma}_{smo} + \vm{\sigma}_{sing}.
\end{equation}

\begin{figure}[!htb]
	\centering
	\includegraphics{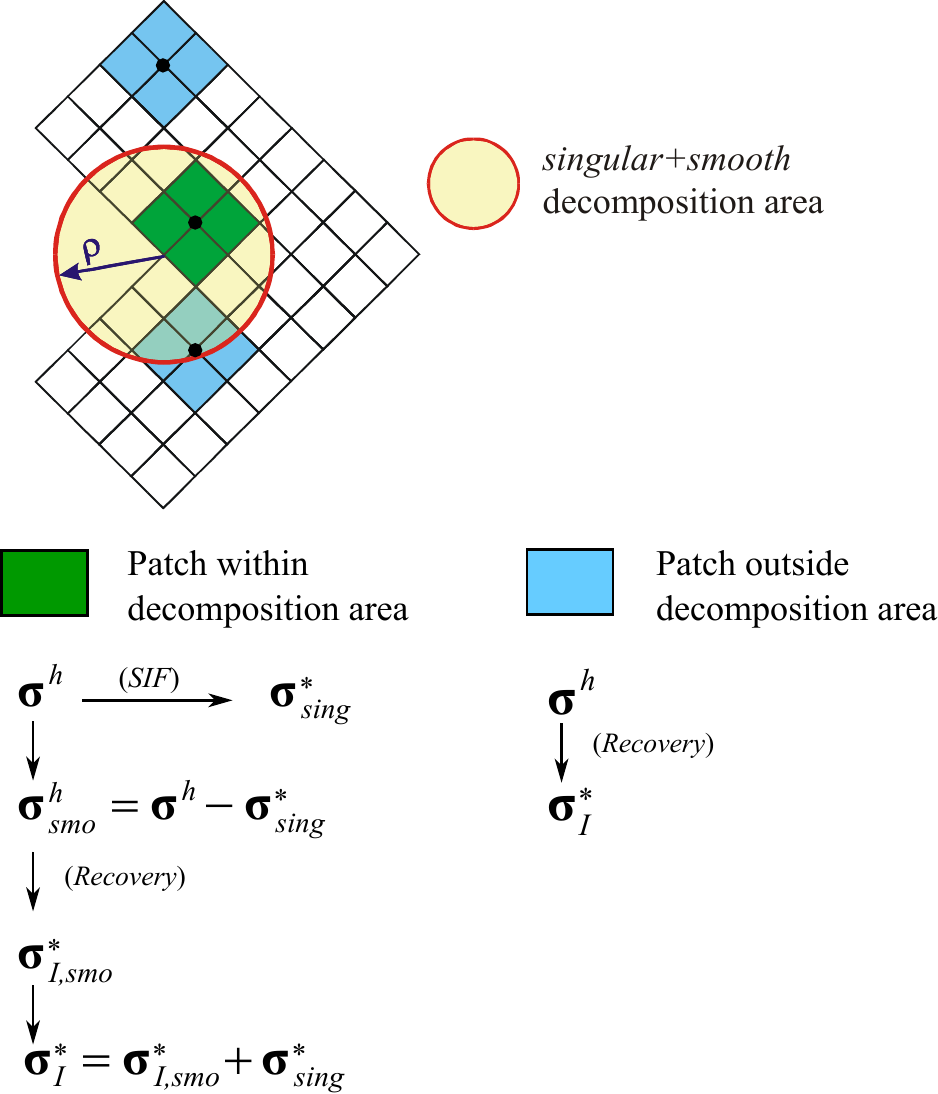}
	\caption{Evaluation of the recovered field at different patches.}
	\label{fig:splitting}
\end{figure}

Then, the recovered field $\vm{\sigma}^*$ required to compute the error estimate given in \Eref{Eq:ZZestimator} can be expressed as the contribution of two recovered stress fields, that is, smooth $\vm{\sigma}^{*}_{smo}$ and singular $\vm{\sigma}^{*}_{sing}$ (see Figure \ref{fig:splitting}):
\begin{equation}  
\vm{\sigma}^* =   \vm{\sigma}^{*}_{smo} + \vm{\sigma}^{*}_{sing}.
\end{equation}

For the recovery of the singular part, the expressions in \Eref{Eq:VnotchStress} which describe the asymptotic fields in the vicinity of the singular point are used. The radius $\rho$ that defines the splitting has to be large enough as to ensure that in the first mesh at least the patch of elements corresponding to the node located at the singular point is within the \emph{singular} + \emph{smooth} decomposition area. To evaluate $\vm{\sigma}^{*}_{sing}$ from  (\ref{Eq:VnotchStress}) we first obtain estimated values of the stress intensity factors $K_{\rm I}$ and $K_{\rm II}$ using the interaction integral as indicated in \cite{rodenasgonzalez2008,rodenasgonzalez2010, yauwang1980}. The recovered part $\vm{\sigma}^{*}_{sing}$ is an equilibrated field as it satisfies the internal equilibrium equations.

Once the field $\vm{\sigma}^{*}_{sing}$ has been evaluated, an FE approximation to the smooth part $\vm{\sigma}^{h}_{smo}$ can be obtained subtracting $\vm{\sigma}^{*}_{sing}$ from the raw FE field:
\begin{equation}   
\vm{\sigma}^{h}_{smo} =   \vm{\sigma}^{h} - \vm{\sigma}^{*}_{sing}.
\end{equation}

Then, the field $\vm{\sigma}^{*}_{smo}$ is evaluated by applying an SPR-C recovery procedure over the field $\vm{\sigma}^{h}_{smo}$.

In order to obtain an equilibrated recovered stress field, the SPR-CX enforces the fulfilment of the equilibrium equations locally on each patch. The constraint equations are introduced via Lagrange multipliers into the linear system used to solve for the coefficients of the polynomial expansion of the recovered stresses on each patch. These include the satisfaction of the:

\begin{itemize}
  \item Internal equilibrium equations.
  \item Boundary equilibrium equation: A point collocation approach is used to impose the satisfaction of a second order approximation to the tractions along the Neumann boundary.  
  \item Compatibility equation: This additional constraint is also imposed to further increase the accuracy of the recovered stress field.  
\end{itemize}

To evaluate the recovered field, quadratic polynomials are used in the patches along the boundary and crack faces, and linear polynomials for the remaining patches. As more information about the solution is available along the boundary, polynomials one degree higher are useful to improve the quality of the recovered stress field.

The enforcement of equilibrium equations provides an equilibrated recovered stress field locally on patches. However, the process used to obtain a continuous field $\vm{\sigma}^*$ shown in (\ref{Eq:conjoint_polynomials}) introduces a small lack of equilibrium as explained in \cite{rodenasgonzalez2010, diezrodenas2007}. The reader is referred to \cite{rodenasgonzalez2010, diezrodenas2007, rodenasgonzalez2010a, rodenasgonzalez2007} for details.

In order to estimate the error in SFEM approximations we can follow a similar procedure. To build the patches we use the topological information of the SFEM discretization. The recovered stress field is evaluated at the centre of the subcells and then projected to the sampling points as explained before. 

\textbf{Remark:} The recovery method proposed in this paper is general, and could also be applied, although this is beyond the scope of this paper, to problems with corner singularities at triple junctions in polycrystalline materials made up of orthotropic grains to estimate the error of extended finite element formulations such as those recently proposed in \cite{menkbordas2010,menkbordas2011b}.

\section{Numerical results}
\label{sec:Results}

In this section, numerical tests considering 2D benchmark problems with exact solution have been used to investigate the quality of the proposed technique. The performance of the technique has been evaluated using the effectivity index of the error in energy norm, both at global and local levels. Globally, we have considered the value of the effectivity index $\theta$ given by 
\begin{equation} \label{Eq:Effectivity}  
\theta  =\frac{\left\| \vm{e}_{es}^{} \right\| }{\left\| \vm{e} \right\| }   
\end{equation}  

where $\|\vm{e}\| $ denotes the exact error in the energy norm, and $\left\| \vm{e}_{es}^{} \right\| $ represents the evaluated error estimate. At element level, the distribution of the local effectivity index $D$, its mean value $m(|D|)$ and standard deviation $\sigma(D)$ have been analysed, as described in \cite{rodenastur2007}: 
\begin{equation} \label{Eq:LocalEffectivity}  
\begin{array}{ccc} 
{D=\theta ^{e} -1} & {\rm if} & {\theta ^{e} \ge 1} \\ 
{D=1-\dfrac{1}{\theta ^{e} }} & {\rm if} & {\theta ^{e} < 1} 
\end{array}
\qquad \qquad {\rm with} \qquad 
\theta ^{e} =\dfrac{\left\| \vm{e}_{es}^{e} \right\| }{\left\| \vm{e}^{e} \right\| }   
\end{equation}  

Note that $\theta^{e} \in (0,1)$ when the error is underestimated and $\theta^{e} \in (1,+\infty)$ when it is overestimated. The definition of $D$ fairly compares the underestimation of the error ($D < 0$) and the overestimation ($D > 0$). The good local behaviour of the estimates results in values of $D$ close to zero. The global effectivity index $\theta$ is used to evaluate global results. The mean value $m(|D|)$ and the standard deviation $\sigma(D)$ of the local effectivity are also used to evaluate the global quality of the error estimator as these parameters are useful to take into account error compensations in the evaluation of $\theta$.

\subsection{ Thick-wall cylinder subjected to an internal pressure.} 
 
\begin{figure}[!htb]
	\centering
	\includegraphics[scale=0.8]{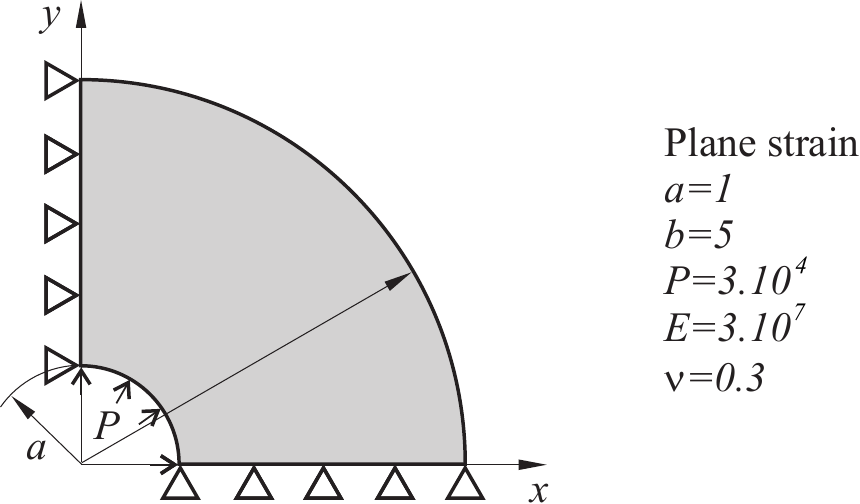}
	\caption{Thick-wall cylinder subjected to internal pressure.}
	\label{fig:cylinder}
\end{figure}

The geometrical model for this problem is shown in Figure~\ref{fig:cylinder}. Due to symmetry conditions, only one part of the section is modelled. Plane strain conditions are assumed.

The exact solution to this problem is given by the following expressions, where for a point $(x,y)$, $c=b/a$, $r=\sqrt{x^{2} +y^{2} }$ and $\phi =\arctan (y/x)$ the radial displacement is given by 
\begin{equation} \label{Eq:uCylinder}  
u_{r} =\frac{P(1+\nu )}{E(c^{2} -1)} \left(r\left(1-2\nu \right) + \frac{b^{2} }{r} \right)  
\end{equation} 
 
Stresses in cylindrical and cartesian coordinates are
\begin{equation} \label{Eq:stressCylinder}  
\begin{array}{cc} 
{\begin{array}{c} {\sigma _{r} =\dfrac{P}{c^{2} -1} \left(1-\dfrac{b^{2} }{r^{2} } \right)} \\ {\sigma _{t} =\dfrac{P}{c^{2} -1} \left(1+\dfrac{b^{2} }{r^{2} } \right)} \end{array}} & 
{\begin{array}{l} {\sigma _{xx} =\sigma _{r} \cos ^{2} (\phi )+\sigma _{t} \sin ^{2} (\phi )} \\ {\sigma _{yy} =\sigma _{r} \sin ^{2} (\phi )+\sigma _{t} \cos ^{2} (\phi )} \\ {\sigma _{xy} =\left(\sigma _{r} -\sigma _{t} \right)\sin (\phi )\cos (\phi )}\\
{\sigma _{z}  =2 \nu \dfrac{P}{c^2-1}  }
\end{array}} \end{array}  
\end{equation} 

\begin{figure}[!htb]
	\centering
	\includegraphics[width=\textwidth]{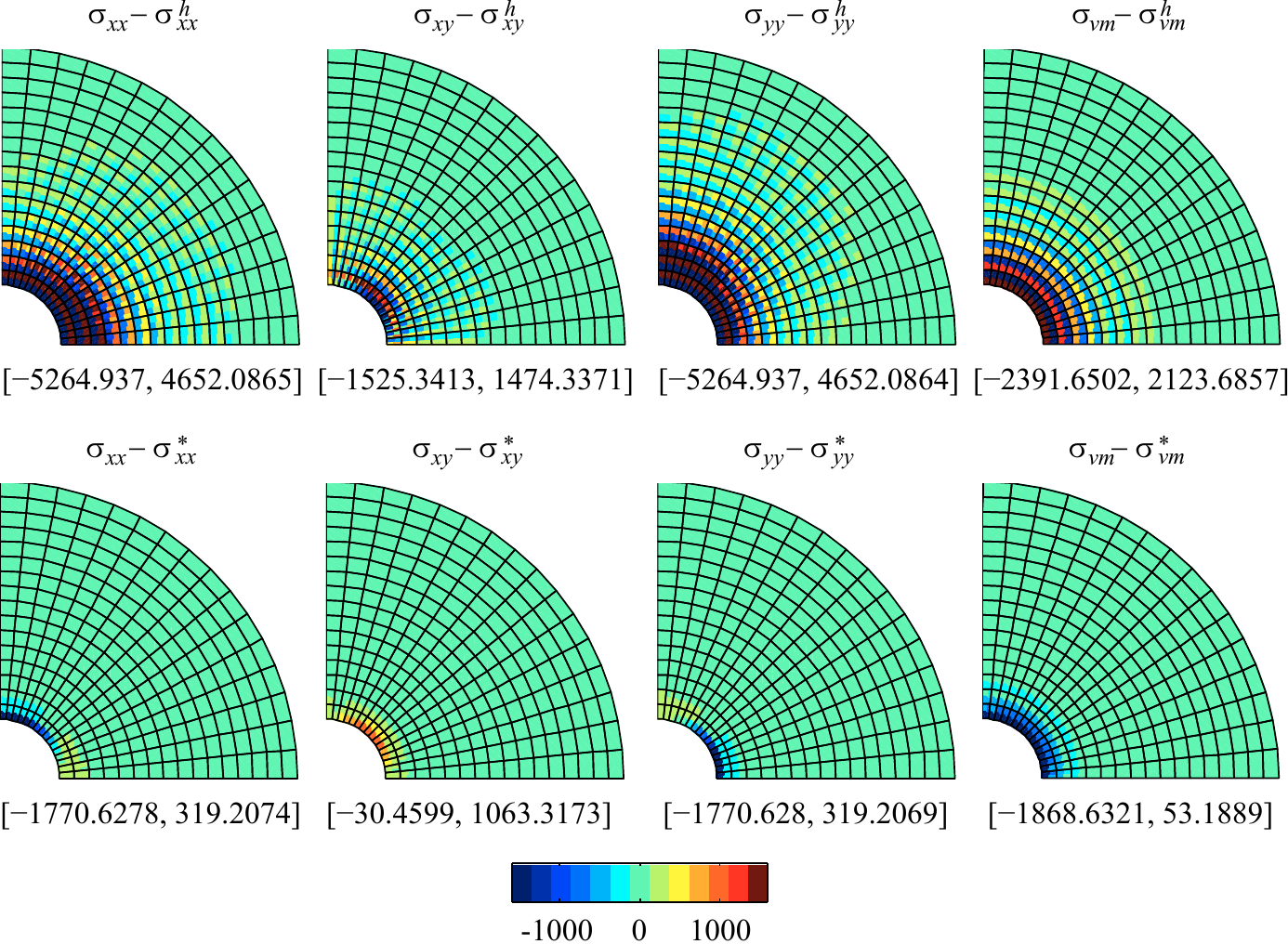}
	\caption{Cylinder under internal pressure. Exact stress error for the SFEM and recovered fields considering the three stress components and the von Mises stress.}
	\label{fig:CylExError}
\end{figure}

A sequence of uniformly refined meshes of linear quadrilateral elements have been used for the analyses. The material parameters are Young's modulus $E=3\times10^7$ and Poisson's ratio $\nu=0.3$. In the case of the SFEM approximation, the element is divided into 4 subcells. However, the influence of the number of subcells on the global/local error level is also studied in a later analysis. Figure \ref{fig:CylExError} shows the exact error for the raw SFEM and the recovered stress fields for the three stress components and the von Mises stress. It can be seen that the error in the recovered field is significantly smaller than the error for the raw stress solution.

\begin{figure}[!htb]
	\centering
	\includegraphics[scale=0.64]{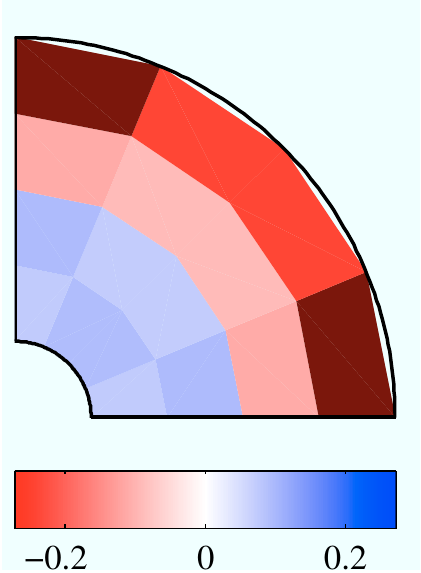}
	\includegraphics[scale=0.64]{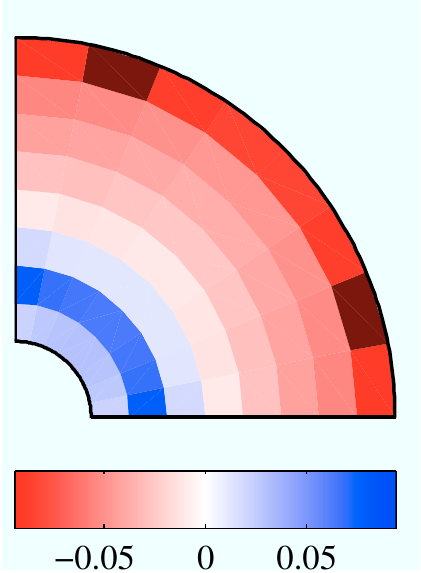}
	\includegraphics[scale=0.64]{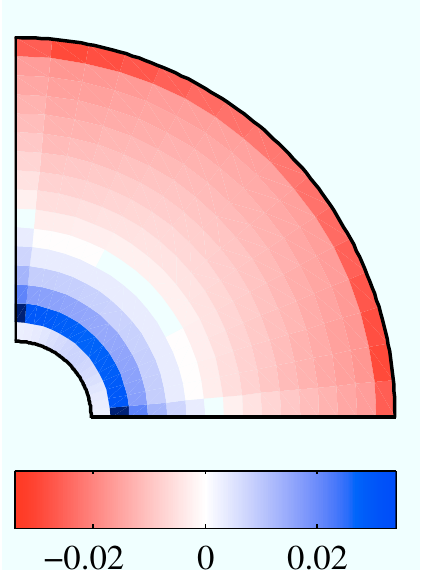}
	\includegraphics[scale=0.64]{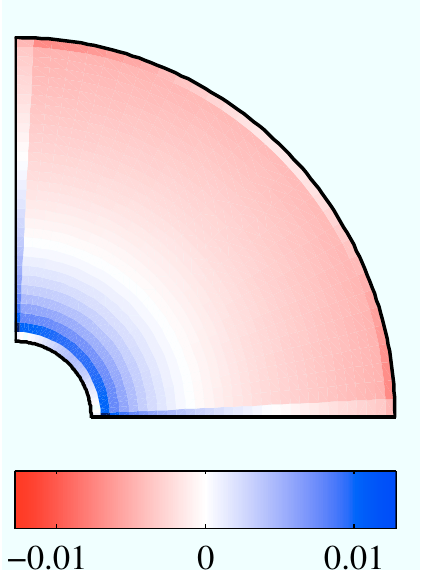}
	\caption{Cylinder under internal pressure. Distribution of the effectivity index $D$.}
	\label{fig:CylD}
\end{figure}

Figure \ref{fig:CylD} shows the distribution of the local effectivity index for a sequence of uniformly refined meshes. The local effectivity values are within a very narrow range and they improve with mesh refinement. The distribution of the $D$ is homogeneous and good results are obtained along the boundary.

\begin{figure}[!htb]
\centering
    \includegraphics{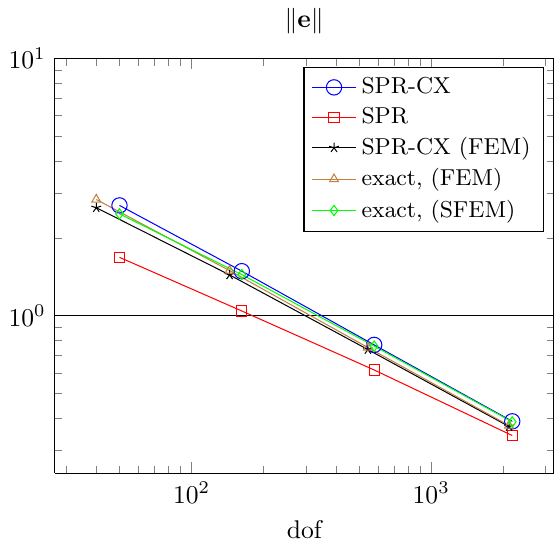}
  \caption{Cylinder under internal pressure. Convergence of the estimated error $\norm{\vm{e}_{es}}$ for the SFEM using 4 subcells (SPR-CX, $s=0.49$), the SFEM without enforcing equilibrium (SPR, $s=0.42$) and the FEM (SPR-CX (FEM), $s=0.49$). The exact error is shown for comparison.}
  \label{fig:CylErrorNorm}
\end{figure}

Figure \ref{fig:CylErrorNorm} shows the convergence of the estimated error in energy norm using two different configurations of the recovery procedure: the curve SPR-CX for an equilibrated recovered field and the curve SPR for a non--equilibrated recovery, the exact error (with a convergence rate $s=0.5$, which equals the optimal convergence rate in the energy norm with respect to the number of degrees of freedom for a smooth problem) is shown for comparison. We have solved the same problem using a standard FEM approximation and an equilibrated recovery technique (curve SPR-CX (FEM)) to estimate the error in that solution, the results are also included in the figure. The FEM values using the SPR-CX recovery converge with a rate of 0.49, while for the SFEM using SPR-CX and standard SPR show an average convergence rate of 0.49 and 0.42, respectively. The equilibrated SFEM and FEM error estimates obtained using the SPR-CX technique are both very close to their corresponding exact errors and converge with the same rate. If we do not consider equilibrium constraints during the recovery (SPR curve) the error is underestimated and the convergence rate is lower. These results clearly show the importance of the use of equilibrated recovery techniques for accurate error estimations.

\begin{figure}[!htb]
\centering
  \includegraphics{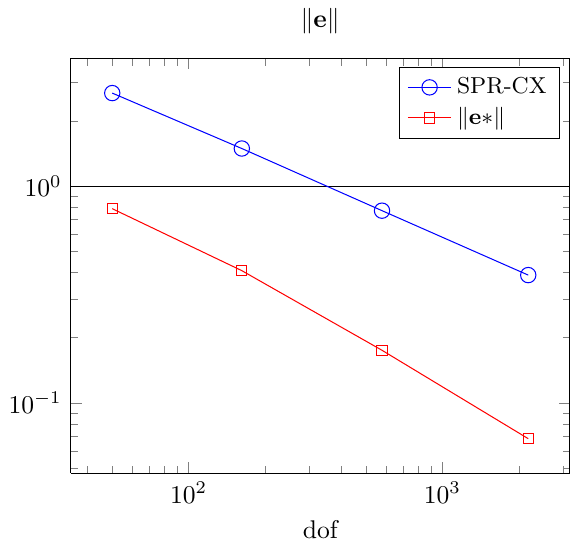}
  \includegraphics{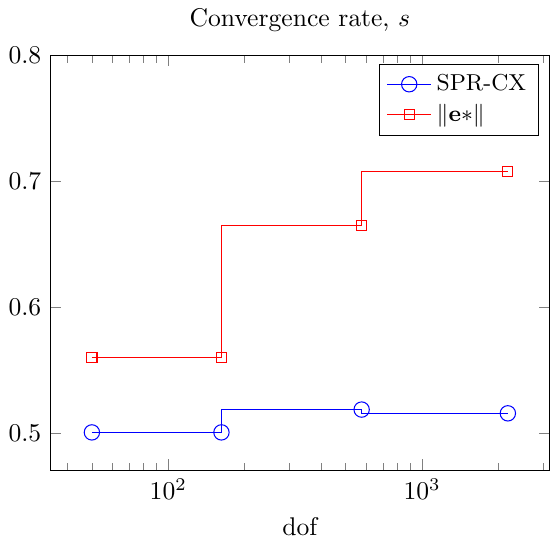}
   \caption{Cylinder under internal pressure. Evolution of the estimated error $\norm{\vm{e}_{es}}$ and convergence rate $s$ for the SFEM using 4 subcells (SPR-CX, $s_{avg}=0.49$) and for the error in the recovered solution ($\norm{\vm{e}^*}$, $s_{avg}=0.65$).}
  \label{fig:CylErrorNormSt}
\end{figure}

In Figure \ref{fig:CylErrorNormSt} we represent the evolution of the estimated error for the SFEM solution using the enhanced recovery (curve SPR-CX) and the exact error of the recovered field $\norm{\vm{e}^*}=\norm{\vm{u}-\vm{u}^*}$. The error in the recovered field has a higher convergence rate, which is in agreement with the expected quality for the field $\vm{\sigma}^*$. According to \cite{zienkiewiczzhu1992a}, this also serves to verify the asymptotic exactness of the proposed error estimator. 

\begin{figure}[!htb]
\centering
\includegraphics{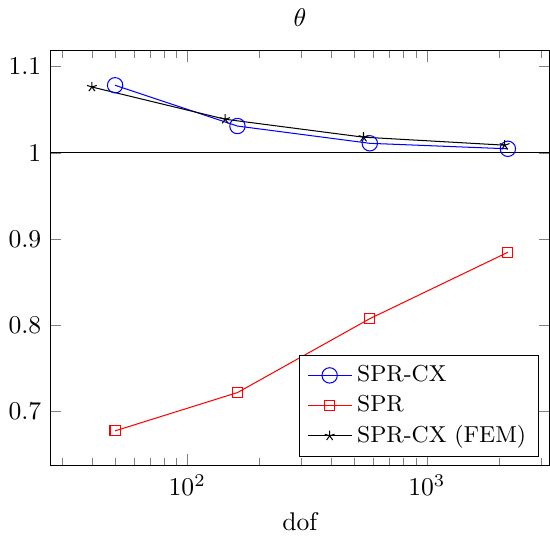}\\
\includegraphics{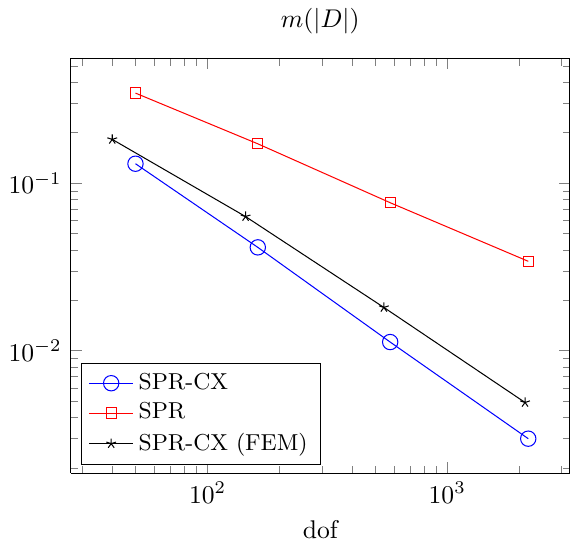}
\includegraphics{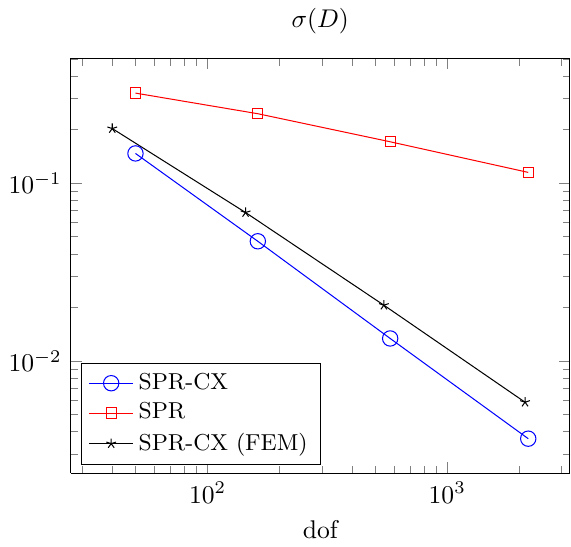}
  \caption{Cylinder under internal pressure. Global indicators $\theta$, $m(|D|)$ and $\sigma(D)$ for SPR-CX, SPR and SPR-CX (FEM).}
  \label{fig:CylGlobal}
\end{figure}

Figure \ref{fig:CylGlobal} shows the evolution of global indicators $\theta$, $m(|D|)$ and $\sigma(D)$ for the SFEM using equilibrated recovery (SPR-CX curve), without equilibrium constraints (SPR curve) and the standard FEM with equilibrium (SPR-CX (FEM) curve). The equilibrated SFEM and the FEM recoveries exhibit similar results, with good effectivity of the error estimator and decrease of  $m(|D|)$ and $\sigma(D)$ for finer meshes. The non--equilibrated SFEM recovery (SPR curve) shows worse results and converges with a slower rate.

\begin{figure}[!htb]
\centering
  \includegraphics{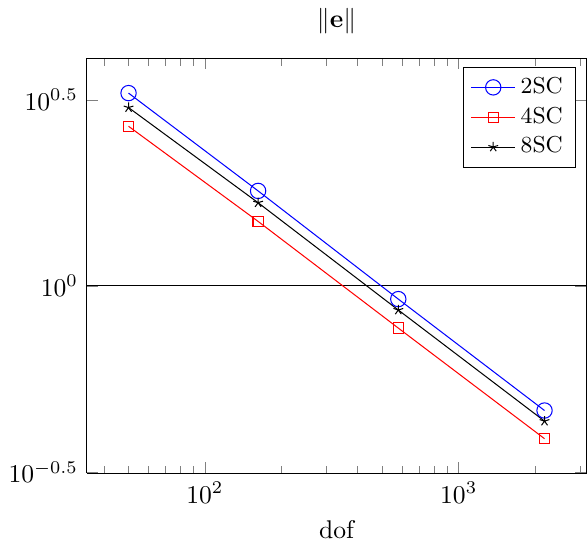}
  \caption{Cylinder under internal pressure. Convergence of the estimated error $\norm{\vm{e}_{es}}$ for elements with two (2SC), four (4SC) and eight (8SC)  subcells.}
  \label{fig:CylErrorNormSC}
\end{figure}

The use of different numbers of subcells for the SFEM approximation is also considered for comparison. Figure \ref{fig:CylErrorNormSC} shows the convergence of the estimated error in energy norm for two, four and eight subcells. All the curves exhibit the same convergence rate $(s=0.49)$, close to the theoretical value $s=0.5$

\begin{figure}[!htb]
\centering
\includegraphics{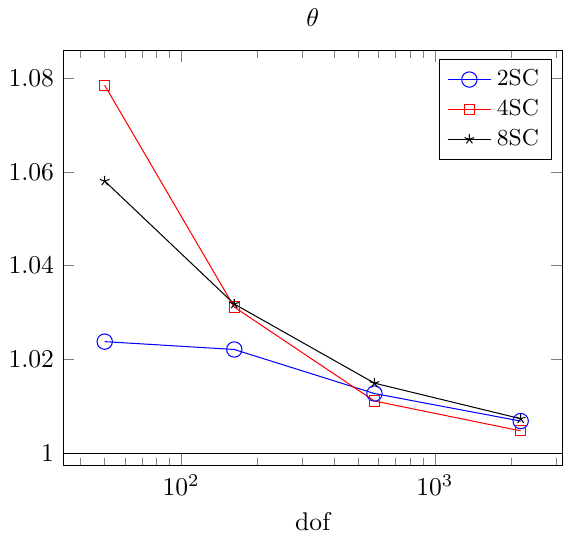}\\
\includegraphics{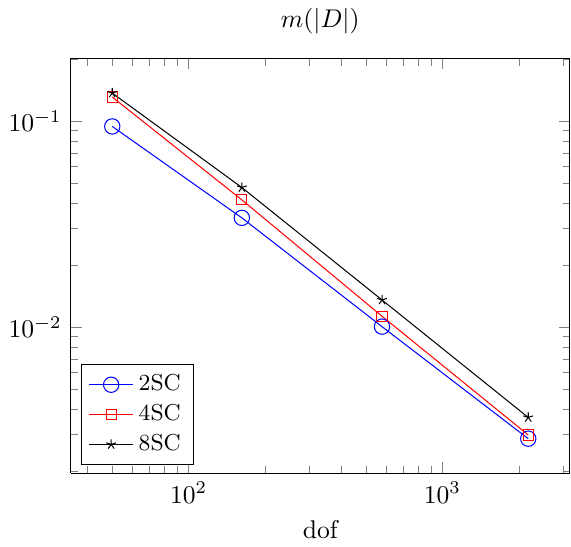}
\includegraphics{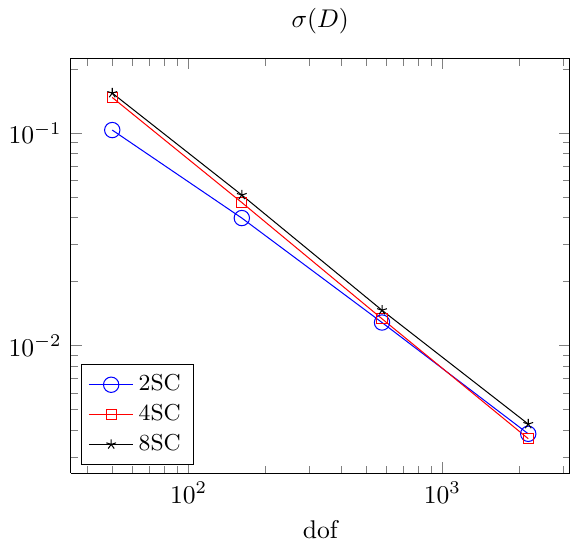}
  \caption{Cylinder under internal pressure. Global indicators $\theta$, $m(|D|)$ and $\sigma(D)$ for elements with two (2SC), four (4SC) and eight (8SC) subcells.}
  \label{fig:CylGlobalSC}
\end{figure}

Figure \ref{fig:CylGlobalSC} shows the evolution of global indicators $\theta$, $m(|D|)$ and $\sigma(D)$ for two, four and eight subcells. The effectivity indices for all the subcells types shown converge asymptotically to the theoretical value and are very accurate ($1.08 > \theta > 1$). The local effectivity index goes to zero at the same rate as shown in the curves $m(|D|)$ and $\sigma(D)$. 

\begin{figure}[!htb]
\centering
\includegraphics{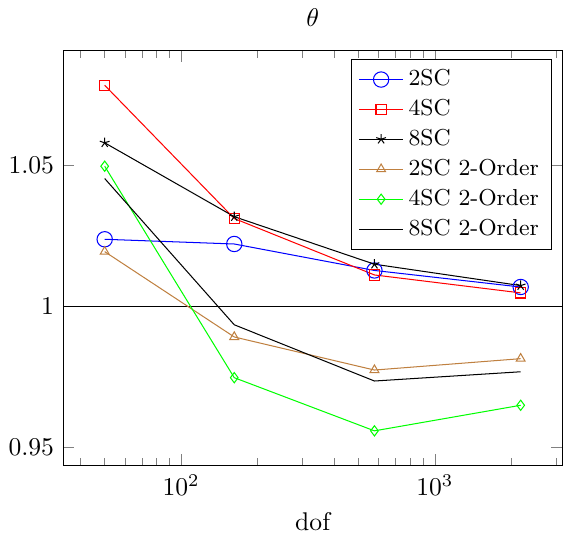}\\
\includegraphics{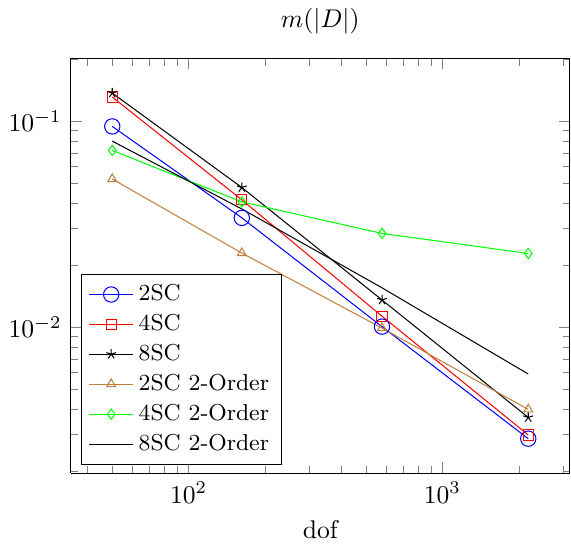}
\includegraphics{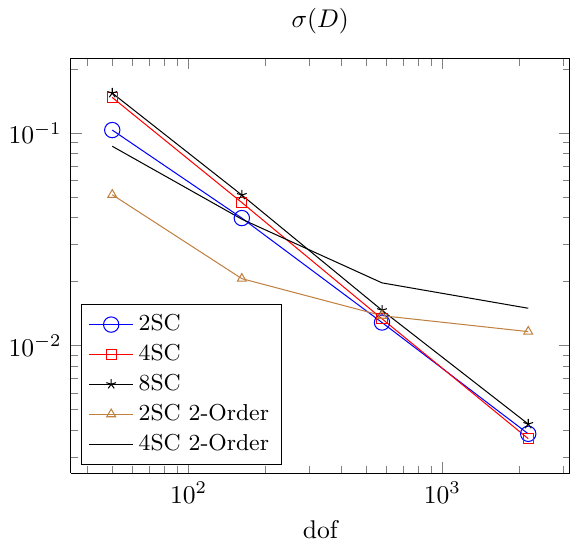}
  \caption{Cylinder under internal pressure.  Global indicators $\theta$, $m(|D|)$ and $\sigma(D)$ for elements with two (2SC), four (4SC) and eight (8SC) subcells using 1st and 2nd order polynomials for the recovery on patches.}
  \label{fig:CylPoly}
\end{figure}

In Figure \ref{fig:CylPoly} we show the influence of the order of the polynomial expansion used for the local recovery on patches. We compare the evolution of the global parameters for first order polynomials,  previously represented in Figure \ref{fig:CylGlobalSC}, with the corresponding curves considering second order polynomials. We can see that the increase of the polynomial order does not produce and improvement of the effectivity. Local behaviour in $m(|D|)$ and $\sigma(D)$  indicate even worse results as we increase the number of degrees of freedom. This is in correspondence with previous results observed in the FEM context \cite{rodenas2001}, where an increase of the polynomial order not necessarily derived in better effectivities.

\subsection{L-shape domain under mode I load.} 
 
\begin{figure}[!htb]
	\centering
	\includegraphics{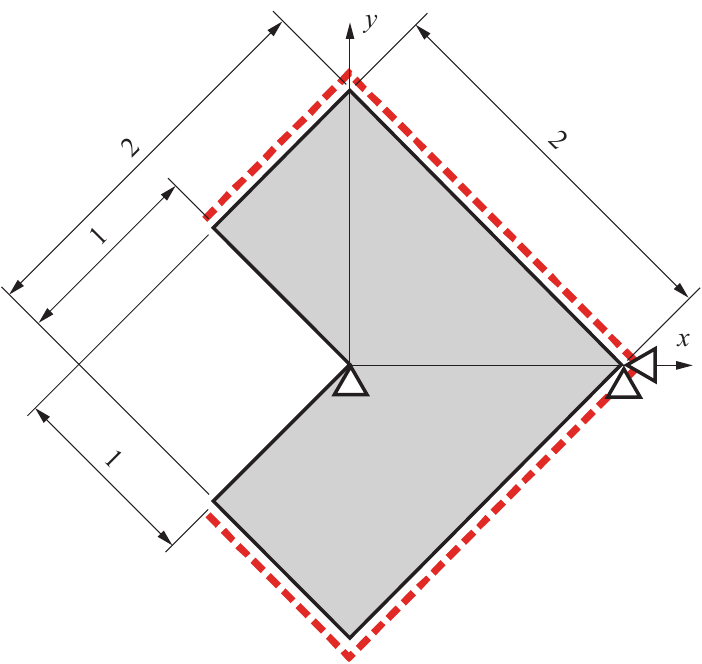}
	\caption{L-shaped domain under mode I load.}
	\label{fig:LShape}
\end{figure}

This singular problem consists of a portion of an infinite L-shaped domain. The model is loaded on the boundary with the tractions corresponding to the first symmetric term of the asymptotic expansion that describes the exact solution under mode I loading conditions around the singular vertex, see Figure~\ref{fig:LShape}. The exact values of boundary tractions on the emphasized boundaries have been imposed in the FE analyses. 

The exact displacement and stress fields for this problem are given by (\ref{Eq:VnotchDisp}, \ref{Eq:VnotchStress}). For $\alpha=3\pi/2$ one obtains $\lambda_{\rm I} =0.544483736782464$, $ \lambda_{\rm II} = 0.908529189846099$ and ${Q}=0.543075578836737$. Exact values of the GSIFs have been taken as $K_{\rm I}=1$ and $ K_{\rm II}=0$. The material parameters are Young's modulus $E~=~1000$, and Poisson's ratio $\nu~=~0.3$. The splitting radius is $\rho=0.5$. To evaluate the SIF we use the equivalent interaction integral defined with a plateau function with radius $0.9$ as indicated in \cite{rodenasgonzalez2008}. As the analytical solution of this problem is singular at the reentrant corner of the plate, we apply the singular+smooth decomposition of the stress field as explained in Section~\ref{section:sprcx}. We use the expression in (\ref{Eq:VnotchStress}) and the estimated GSIF evaluated using the interaction integral to obtain an approximation of the singular part $\sigma_{sing}$.

\begin{figure}[!htb]
	\centering
	\includegraphics[scale=0.64]{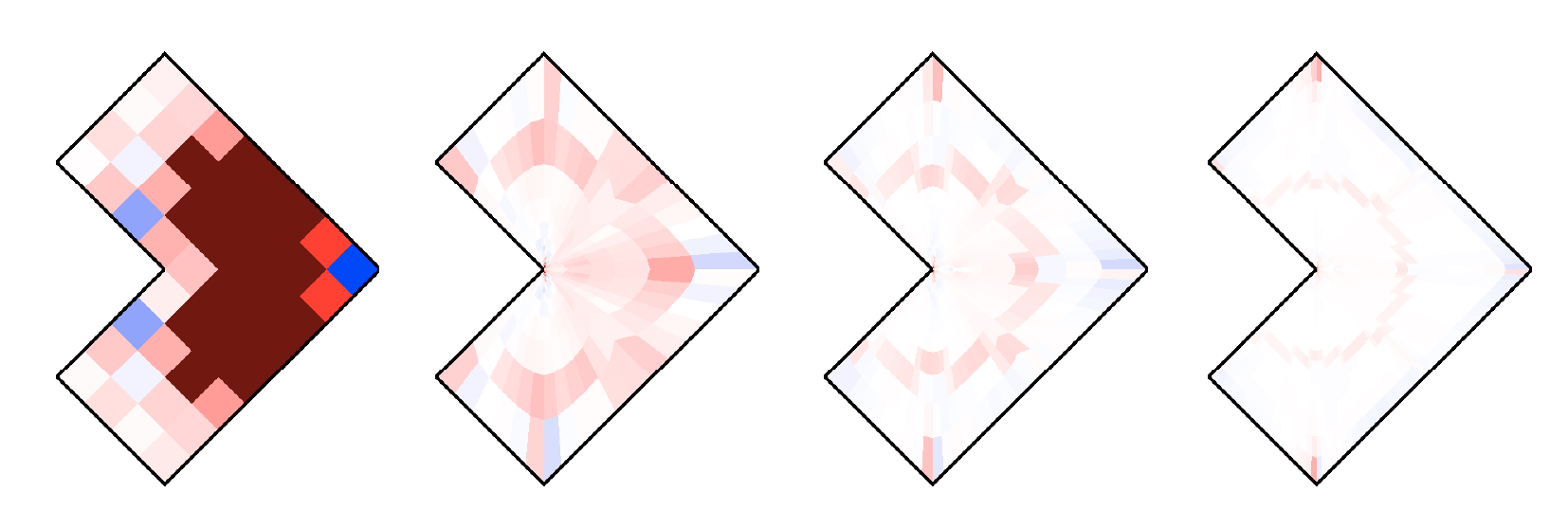}
	\vspace{10pt}
	\includegraphics[scale=0.7]{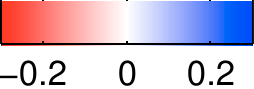}\\
	
	\caption{L-shaped domain under mode I load. Distribution of the effectivity index $D$.}
	\label{fig:LSHPD}
\end{figure}

\begin{figure}[!htb]
	\centering
	\includegraphics[scale=0.64]{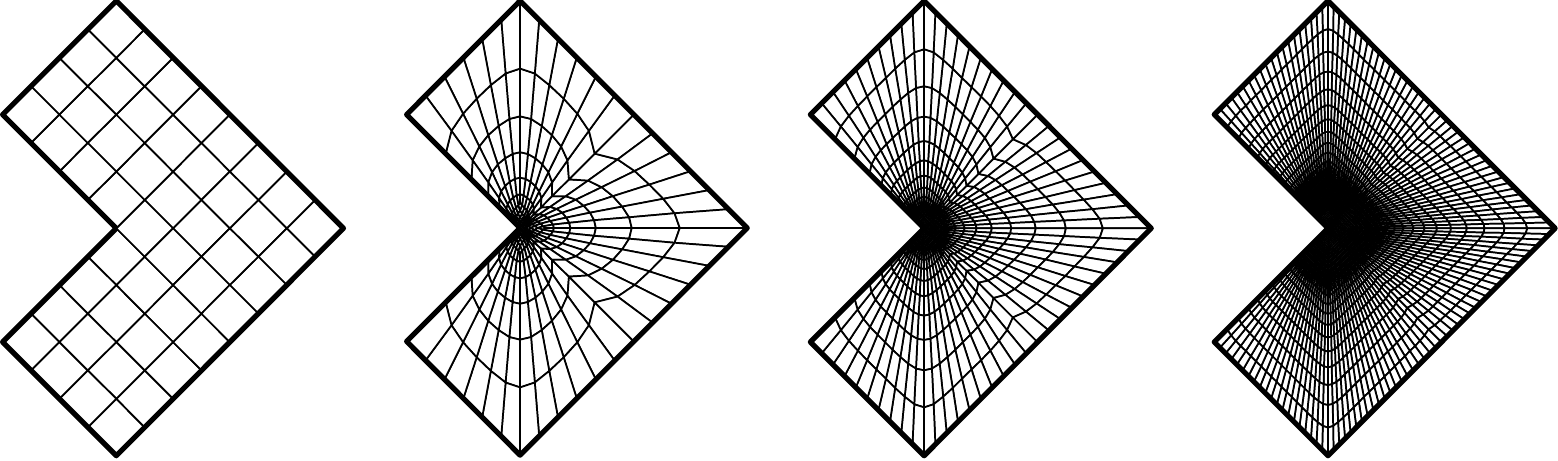}
	\caption{L-shaped domain under mode I load. Sequence of graded meshes.}
	\label{fig:LSHPmesh}
\end{figure}

Figure \ref{fig:LSHPD} shows the distribution of the local effectivity index for the sequence of graded meshes in Figure \ref{fig:LSHPmesh}. In the figure, the local index $D$ decreases with the refinement of the meshes and the obtained values are within a narrow range. The first mesh of the sequence is an uniform mesh. This kind of meshes is known to produce pollution error for singular problems. The highly underestimated areas at the right of the plate are explained by this pollution error, which cannot be controlled with local techniques.

\begin{figure}[!htb]
\centering
\includegraphics{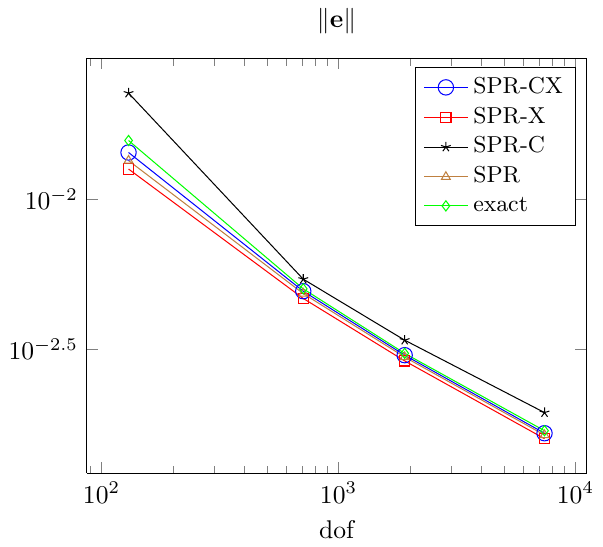}
%
  \caption{L-shaped domain under mode I load. Convergence of the estimated error $\norm{\vm{e}_{es}}$ for for different configurations of the recovery technique: SPR-CX, SPR-X, SPR-C, SPR. The exact error is shown for comparison.}
  \label{fig:LSHPErrorNorm}
\end{figure}

Figure \ref{fig:LSHPErrorNorm} shows the convergence of the estimated error in energy norm for different configurations of the recovery procedure: SPR-CX that consider equilibrium and stress decomposition, SPR-X that considers stress decomposition only, SPR-C that consider equilibrium only, and a conventional SPR. The exact error is shown for comparison. 

\begin{figure}[!htb]
\centering
\includegraphics{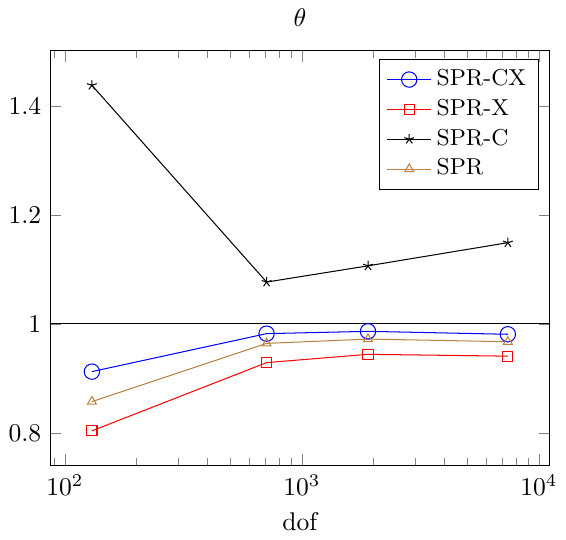}\\

 \begin{minipage}[c] {0.49\textwidth}
 \includegraphics{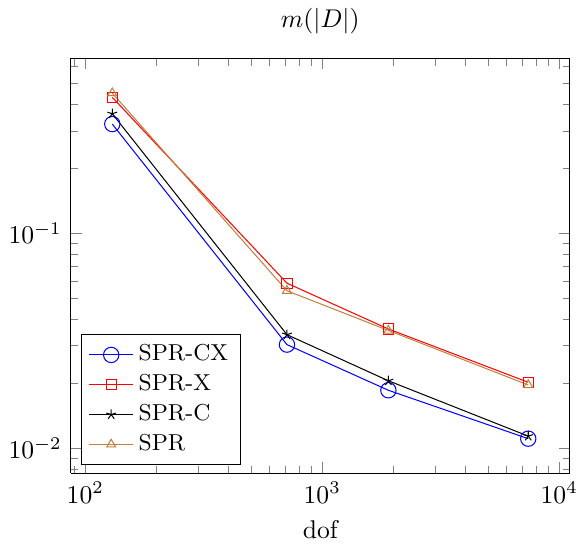}
 \end{minipage} 
 \begin{minipage}[c]{0.49\textwidth}
   \includegraphics{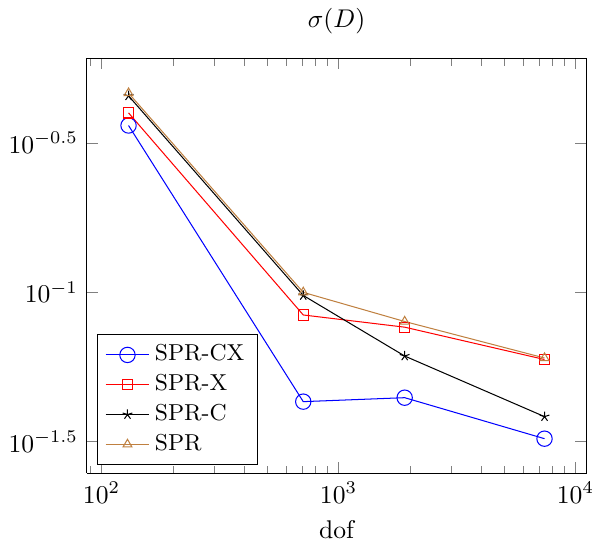}
 \end{minipage}

  \caption{L-shaped domain under mode I load. Global indicators $\theta$, $m(|D|)$ and $\sigma(D)$ for different configurations of the recovery technique: SPR-CX, SPR-X, SPR-C, SPR. }
  \label{fig:LSHPGlobal}
\end{figure}

\begin{figure}[!htb]
	\centering
        \begin{minipage}[c]{0.49\textwidth}
         \centering
	  SPR-CX\\
	  \includegraphics{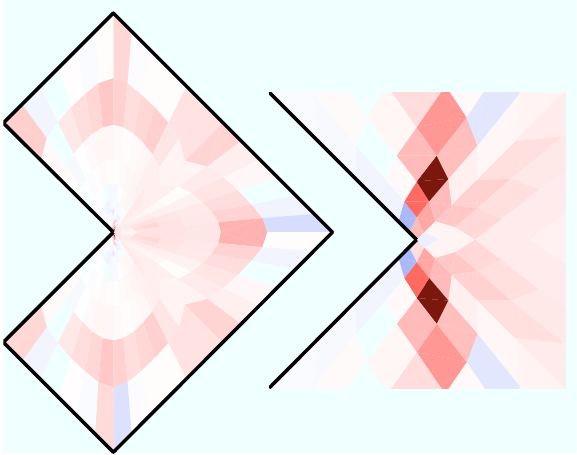}
        \end{minipage}
        \begin{minipage}[c]{0.49\textwidth}
         \centering
	  SPR-X\\
	  \includegraphics{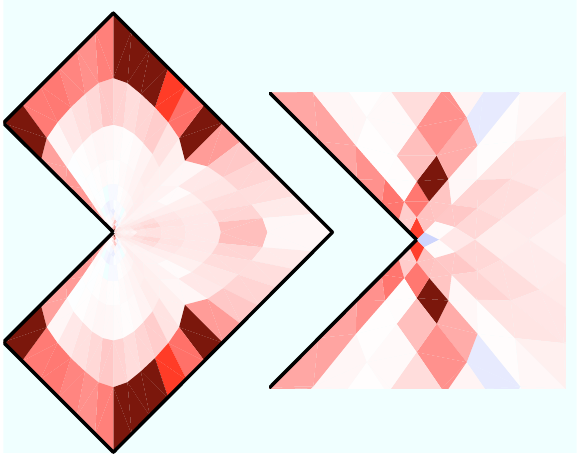}
        \end{minipage} \\
	\vspace{20pt}
	\begin{minipage}[c]{0.49\textwidth}
         \centering
	  SPR-C \\
	  \includegraphics{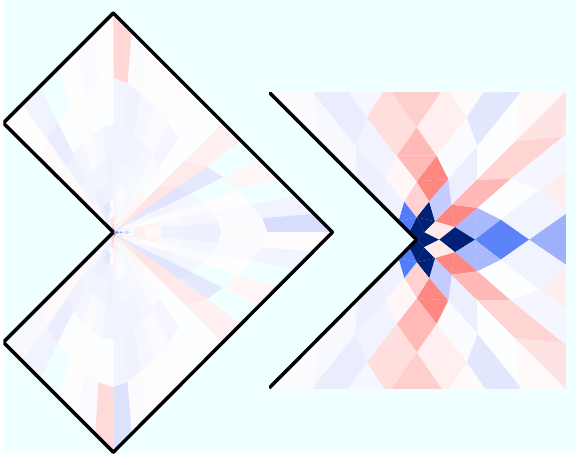}
        \end{minipage}
	\begin{minipage}[c]{0.49\textwidth}
         \centering
	  SPR\\
	  \includegraphics{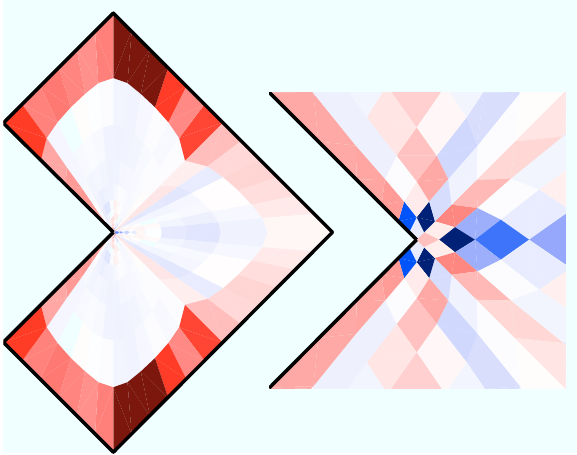}
        \end{minipage}

	\vspace{20pt}
	\includegraphics[scale=0.7]{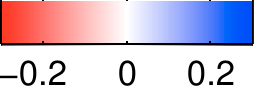}\\
	
	\caption{L-shaped domain under mode I load. Distribution of the effectivity index $D$ for the different configurations of the recovery technique. SPR-CX: equilibrium and stress decomposition, SPR-X: stress decomposition, SPR-C: equilibrium, and a conventional SPR.}
	\label{fig:LSHPDzoom}
\end{figure}

Figure \ref{fig:LSHPGlobal} shows the evolution of global indicators $\theta$, $m(|D|)$ and $\sigma(D)$ for the different configurations: SPR-CX, SPR-X, SPR-C, SPR. The best results are for the SPR-CX. The SPR-C and the SPR cannot accurately recover the field close to the singularity. The SPR seems to provide good global effectivity results, however, this is only due to compensation between underestimated and overestimated areas. The real behaviour is clear when analysing the evolution of $m(|D|)$ and $\sigma(D)$ and even more patent if we represent $D$ as seen in Figure \ref{fig:LSHPDzoom}, where we represent the results for the different configurations of the recovery technique. The effect of equilibrium enforcement and singular decomposition is clearly shown. A more homogeneous distribution of $D$ is obtained when using SPR-CX. SPR and SPR-X are not equilibrated along the boundary and therefore they underestimate the error in the elements along it. SPR and SPR-C overestimates the error in the vicinity of the reentrant corner.

\section{Conclusions}
\label{sec:Conclusions}

In this paper, an \emph{a posteriori} recovery-based error estimator which makes use of a modified version of the SPR technique previously used in the FEM and XFEM contexts has been adapted to the smoothed finite element method. The recovery technique considers the local enforcement of equilibrium equations and the \emph{singular} + \emph{smooth} decomposition of the stress field for singular problems. The technique has been applied to the cell-based smoothed FEM  but could also be used with the node-based and edge-based smoothed FEM implementations.

The numerical results presented in this paper show that the method yields accurate estimations of the error in the energy norm both locally and globally. It can be inferred that enforcing equilibrium constraint is required to obtain accurate results. The influence in the recovery of the number of subcells used to formulate the smoothed finite elements is also studied, showing that the technique performs adequately for the different configurations. The error estimator based on the use of the SPR-CX recovery technique accurately captures the discretization error both in smooth and singular problems. Moreover, it could be used to guide \emph{h}-adaptive refinements and the recovered field $\bm{\sigma}^*$ can be used as an enhanced solution, more accurate than the stress field provided by the approximation. 

Future work includes the comparison of the proposed estimators with other recently developed error estimators for the extended finite element method when dealing with singular problems \cite{bordasduflot2007,duflotbordas2008,bordasduflot2008} for three dimensional fracture problems, the focus of our current research. We will also analyse the behaviour of strain smoothing for real-life three dimensional fracture mechanics problems, which is the topic of the EPSRC project which funded this work.

\section{Acknowledgements}

St\'{e}phane Bordas would like to thank the partial financial support of the Royal Academy of Engineering and of the Leverhulme Trust for his Senior Research Fellowship \emph{Towards the next generation surgical simulators} as well as the financial support for Octavio A. Gonz\'{a}lez-Estrada and St{\' e}phane Bordas from the UK Engineering Physical Science Research Council (EPSRC) under grant EP/G042705/1 \emph{Increased Reliability for Industrially Relevant Automatic Crack Growth Simulation with the eXtended Finite Element Method}. St\'{e}phane Bordas also thanks partial financial support of the European Research Council Starting Independent Research Grant (ERC Stg grant agreement No. 279578). \\

This work has been carried out within the framework of the research project DPI2010-20542 of the Ministerio de Ciencia e Innovaci\'{o}n (Spain). The financial support from Universitat Polit\`{e}cnica de Val\`{e}ncia and Generalitat Valenciana are also acknowledged.

\bibliographystyle{wileyj}
\bibliography{library}

\begin{thebibliography}{10}
\providecommand{\url}[1]{\texttt{#1}}
\providecommand{\urlprefix}{URL }
\expandafter\ifx\csname urlstyle\endcsname\relax
  \providecommand{\doi}[1]{doi:\discretionary{}{}{}#1}\else
  \providecommand{\doi}{doi:\discretionary{}{}{}\begingroup
  \urlstyle{rm}\Url}\fi

\bibitem{liudai2006}
Liu GR, Dai KY, Nguyen TT. {A Smoothed Finite Element Method for Mechanics
  Problems}. \emph{Computational Mechanics}  May 2006; \textbf{39}(6):859--877,
  \doi{10.1007/s00466-006-0075-4}.

\bibitem{liunguyen2007}
Liu GR, Nguyen TT, Dai KY, Lam KY. Theoretical aspects of the smoothed finite
  element method ({SFEM}). \emph{International Journal for Numerical Methods in
  Engineering}  2007; \textbf{71}(8):902--930.

\bibitem{nguyen-xuanbordas2008}
Nguyen-Xuan H, Bordas SPA, Nguyen-Dang H. {Smooth finite element methods:
  convergence, accuracy and properties}. \emph{International Journal for
  Numerical Methods in Engineering}  2008; \textbf{74}(2):175--208,
  \doi{10.1002/nme}.

\bibitem{bordasnatarajan2010}
Bordas SPA, Natarajan S. {On the approximation in the smoothed finite element
  method (SFEM)}. \emph{International Journal for Numerical Methods in
  Engineering}  2010; \textbf{81}(5):660--670, \doi{10.1002/nme}.

\bibitem{zhangliu2008}
Zhang HH, Liu SJ, Li LX. On the smoothed finite element method.
  \emph{International Journal for Numerical Methods in Engineering}  2008;
  \textbf{76}(8):1285--1295, \doi{10.1002/nme.2460}.

\bibitem{nguyen-thoiliu2009}
Nguyen-Thoi T, Liu G, Lam K, Zhang G. A face-based smoothed finite element
  method {(FS-FEM)} for 3{D} linear and nonlinear solid mechanics using 4-node
  tetrahedral elements. \emph{International Journal for Numerical Methods in
  Engineering}  2009; \textbf{78}:324--353.

\bibitem{liunguyen-thoi2009b}
Liu G, Nguyen-Thoi T, Lam K. An edge-based smoothed finite element method
  {(ES-FEM)} for static, free and forced vibration analyses of solids.
  \emph{Journal of Sound and Vibration}  2009; \textbf{320}:1100--1130.

\bibitem{liunguyen-thoi2009}
Liu G, Nguyen-Thoi T, Nguyen-Xuan H, Lam K. A node based smoothed finite
  element method {(NS-FEM)} for upper bound solution to solid mechanics
  problems. \emph{Computers and Structures}  2009; \textbf{87}:14--26.

\bibitem{liu2010}
Liu G. \emph{Smoothed Finite Element Methods}. CRC Press, 2010.

\bibitem{liunguyen-xuan2010}
Liu G, Nguyen-Xuan H, Nguyen-Thoi T. A theoretical study on the smoothed {FEM
  (SFEM)} models: {P}roperties, accuracy and convergence rates.
  \emph{International Journal for Numerical Methods in Biomedical Engineering}
  2010; \textbf{84}:1222--1256.

\bibitem{nguyenliu2007}
Nguyen T, Liu G, Dai K, Lam K. Selective smoothed finite element method.
  \emph{Tsinghua {S}cience \& {T}echnology}  2007; \textbf{12}:497--508.

\bibitem{hungbordas2009}
Hung NX, Bordas S, Hung N. Addressing volumetric locking and instabilities by
  selective integration in smoothed finite element. \emph{Communications in
  Numerical Methods in Engineering}  2009; \textbf{25}:19--34.

\bibitem{nguyen-xuanrabczuk2008}
Nguyen-Xuan H, Rabczuk T, Bordas S, Debongnie JF. A smoothed finite element
  method for plate analysis. \emph{Computer Methods in Applied Mechanics and
  Engineering}  2008; \textbf{197}:1184--1203.

\bibitem{nguyenrabczuk2008}
Nguyen NT, Rabczuk T, Nguyen-Xuan H, Bordas S. A smoothed finite element method
  for shell analysis. \emph{Computer Methods in Applied Mechanics and
  Engineering}  2008; \textbf{198}:165--177.

\bibitem{bordasrabczuk2010}
Bordas SPA, Rabczuk T, Hung NX, Nguyen VP, Natarajan S, Bog T, Quan DM, Hiep
  NV. {Strain smoothing in FEM and XFEM}. \emph{Computers \& Structures}  Dec
  2010; \textbf{88}(23-24):1419--1443, \doi{10.1016/j.compstruc.2008.07.006}.

\bibitem{bordasnatarajan2011}
Bordas SP, Natarajan S, Kerfriden P, Augarde CE, Mahapatra DR, Rabczuk T, Pont
  SD. On the performance of strain smoothing for quadratic and enriched finite
  element approximations ({XFEM/GFEM/PUFEM}). \emph{International Journal for
  Numerical Methods in Biomedical Engineering}  2011; \textbf{86}:637--666.

\bibitem{liunguyen-thoi2009a}
Liu G, Nguyen-Thoi T, Nguyen-Xuan H, Dai K, Lam K. On the essence and the
  evaluation of the shape functions for the smoothed finite element method
  ({SFEM}). \emph{International Journal for Numerical Methods in Engineering}
  2009; \textbf{77}:1863--1869, \doi{10.1002/nme.2587}.

\bibitem{stroubouliszhang2006}
Strouboulis T, Zhang L, Wang D, Babu\v{s}ka I. {A posteriori error estimation
  for generalized finite element methods}. \emph{Computer Methods in Applied
  Mechanics and Engineering}  2006; \textbf{195}(9-12):852--879.

\bibitem{bordasduflot2007}
Bordas SPA, Duflot M. {Derivative recovery and a posteriori error estimate for
  extended finite elements}. \emph{Computer Methods in Applied Mechanics and
  Engineering}  2007; \textbf{196}(35-36):3381--3399.

\bibitem{xiaokarihaloo2004}
Xiao QZ, Karihaloo BL. {Statically admissible stress recovery using the moving
  least squares technique}. \emph{Progress in Computational Structures
  Technology}, Topping BHV, Soares CAM (eds.), Saxe-Coburg Publications:
  Stirling, Scotland, 2004; 111--138.

\bibitem{rodenasgonzalez2008}
R\'{o}denas JJ, Gonz\'{a}lez-Estrada OA, Taranc\'{o}n JE, Fuenmayor FJ. {A
  recovery-type error estimator for the extended finite element method based on
  singular+smooth stress field splitting}. \emph{International Journal for
  Numerical Methods in Engineering}  2008; \textbf{76}(4):545--571,
  \doi{10.1002/nme.2313}.

\bibitem{panetierladeveze2010}
Panetier J, Ladev\`{e}ze P, Chamoin L. {Strict and effective bounds in
  goal-oriented error estimation applied to fracture mechanics problems solved
  with XFEM}. \emph{International Journal for Numerical Methods in Engineering}
   2010; \textbf{81}(6):671--700.

\bibitem{nguyen-thoiliu2011}
Nguyen-Thoi T, Liu G, Nguyen-Xuan H, Nguyen-Tran C. {Adaptive analysis using
  the node-based smoothed finite element method (NS-FEM)}. \emph{International
  Journal for Numerical Methods in Biomedical Engineering}  2011;
  \textbf{27}(2):198--218, \doi{10.1002/cnm}.

\bibitem{zienkiewiczzhu1987}
Zienkiewicz OC, Zhu JZ. {A simple error estimator and adaptive procedure for
  practical engineering analysis}. \emph{International Journal for Numerical
  Methods in Engineering}  1987; \textbf{24}(2):337--357.

\bibitem{rodenasgonzalez2010}
R\'{o}denas JJ, Gonz\'{a}lez-Estrada OA, D\'{\i}ez P, Fuenmayor FJ. {Accurate
  recovery-based upper error bounds for the extended finite element framework}.
  \emph{Computer Methods in Applied Mechanics and Engineering}  2010;
  \textbf{199}(37-40):2607--2621.

\bibitem{williams1952}
Williams ML. {Stress singularities resulting from various boundary conditions
  in angular corners of plate in extension}. \emph{Journal of Applied
  Mechanics}  1952; \textbf{19}:526--534.

\bibitem{szabobabuska1991}
Szab\'{o} BA, Babu\v{s}ka I. \emph{{Finite Element Analysis}}. John Wiley \&
  Sons: New York, 1991.

\bibitem{barber2010}
Barber JR. \emph{{Elasticity. Series: Solid Mechanics and its applications}}.
  3rd edn., Springer: Dordrecht, 2010.

\bibitem{chenwu2001}
Chen JS, Wu CT, Yoon S, You Y. A stabilized conforming nodal integration for
  {Galerkin} mesh-free methods. \emph{Int. J. Numer. Meth. Engng.}  2001;
  \textbf{50}:435--466.

\bibitem{yoomoran2004}
Yoo J, Moran B, Chen J. Stabilized conforming nodal integration in the natural
  element method. \emph{International Journal for Numerical Methods in
  Engineering}  2004; \textbf{60}:861--890.

\bibitem{zienkiewiczzhu1992}
Zienkiewicz OC, Zhu JZ. {The superconvergent patch recovery and a posteriori
  error estimates. Part 1: The recovery technique}. \emph{International Journal
  for Numerical Methods in Engineering}  1992; \textbf{33}(7):1331--1364.

\bibitem{zienkiewiczzhu1992a}
Zienkiewicz OC, Zhu JZ. {The superconvergent patch recovery and a posteriori
  error estimates. Part 2: Error estimates and adaptivity}. \emph{International
  Journal for Numerical Methods in Engineering}  1992;
  \textbf{33}(7):1365--1382.

\bibitem{duflotbordas2008}
Duflot M, Bordas SPA. {A posteriori error estimation for extended finite
  elements by an extended global recovery}. \emph{International Journal for
  Numerical Methods in Engineering}  2008; \textbf{76}:1123--1138,
  \doi{10.1002/nme}.

\bibitem{bordasduflot2008}
Bordas SPA, Duflot M, Le P. {A simple error estimator for extended finite
  elements}. \emph{Communications in Numerical Methods in Engineering}  2008;
  \textbf{24}(11):961--971.

\bibitem{rodenastur2007}
R\'{o}denas JJ, Tur M, Fuenmayor FJ, Vercher A. {Improvement of the
  superconvergent patch recovery technique by the use of constraint equations:
  the SPR-C technique}. \emph{International Journal for Numerical Methods in
  Engineering}  2007; \textbf{70}(6):705--727, \doi{10.1002/nme.1903}.

\bibitem{diezrodenas2007}
D\'{\i}ez P, R\'{o}denas JJ, Zienkiewicz OC. {Equilibrated patch recovery error
  estimates: simple and accurate upper bounds of the error}.
  \emph{International Journal for Numerical Methods in Engineering}  2007;
  \textbf{69}(10):2075--2098, \doi{10.1002/nme}.

\bibitem{wibergabdulwahab1993}
Wiberg NE, Abdulwahab F. {Patch recovery based on superconvergent derivatives
  and equilibrium}. \emph{International Journal for Numerical Methods in
  Engineering}  Aug 1993; \textbf{36}(16):2703--2724,
  \doi{10.1002/nme.1620361603}.

\bibitem{blackerbelytschko1994}
Blacker T, Belytschko T. {Superconvergent patch recovery with equilibrium and
  conjoint interpolant enhancements}. \emph{International Journal for Numerical
  Methods in Engineering}  1994; \textbf{37}(3):517--536.

\bibitem{yauwang1980}
Yau J, Wang S, Corten H. {A mixed-mode crack analysis of isotropic solids using
  conservation laws of elasticity}. \emph{Journal of Applied Mechanics}  1980;
  \textbf{47}(2):335--341.

\bibitem{rodenasgonzalez2010a}
R\'{o}denas JJ, Gonz\'{a}lez-Estrada OA, Fuenmayor FJ, Chinesta F. {Upper
  bounds of the error in X-FEM based on a moving least squares (MLS) recovery
  technique}. \emph{9th. World Congress on Computational Mechanics (WCCM9).
  4th. Asian Pacific Congress on Computational Methods (APCOM2010)}, Khalili N,
  Valliappan S, Li Q, Russell A (eds.), Centre for Infrastructure Engineering
  and Safety, 2010.

\bibitem{rodenasgonzalez2007}
R\'{o}denas JJ, Gonz\'{a}lez-Estrada OA, D\'{\i}ez P, Fuenmayor FJ. {Upper
  bounds of the error in the extended finite element method by using an
  equilibrated-stress patch recovery technique}. \emph{International Conference
  on Adaptive Modeling and Simulation. ADMOS 2007}, International Center for
  Numerical Methods in Engineering (CIMNE), 2007; 210--213.

\bibitem{menkbordas2010}
Menk A, Bordas S. Numerically determined enrichment function for the extended
  finite element method and applications to bi-material anisotropic fracture
  and polycrystals. \emph{International Journal for Numerical Methods in
  Engineering}  2010; \textbf{83}:805--828.

\bibitem{menkbordas2011b}
Menk A, Bordas S. Crack growth calculations in solder joints based on
  microstructural phenomena with {X-FEM}. \emph{Computational Materials
  Science}  2011; (3):1145--1156.

\bibitem{rodenas2001}
R\'odenas JJ. {Error de discretizaci\'{o}n en el c\'{a}lculo de sensibilidades
  mediante el m\'etodo de los elementos finitos}. Ph{D} {T}hesis, Universidad
  Polit\'{e}cnica de Valencia 2001.

\end{thebibliography}

\end{document}